\documentclass[11pt,a4paper]{article}

\usepackage{fancyhdr}
\usepackage{amsmath}
\usepackage{bm}
\usepackage{esint}
\usepackage{amsfonts}
\usepackage{amsthm}
\usepackage{amssymb}
\usepackage{amsbsy}
\usepackage{verbatim}
\usepackage{graphicx}
\usepackage{url}
\usepackage{mathrsfs}
\usepackage[hmargin = 1in,vmargin=.75in]{geometry}
\usepackage{color}
\usepackage{amstext}
\usepackage{stmaryrd}
\usepackage{enumerate}
\usepackage{algpseudocode}
\usepackage{algorithm}
\usepackage{pifont}
\usepackage{prettyref}
\usepackage{hyperref}
\hypersetup{
    colorlinks = true,
}

\font\msbm=msbm10

\numberwithin{equation}{section}

\theoremstyle{plain}
\newtheorem{theorem}{Theorem}[section]
\newtheorem{lemma}[theorem]{Lemma}

\newtheorem{proposition}[theorem]{Proposition}

\def\mathbb#1{\hbox{\msbm{#1}}}

\newcommand{\be}{\boldsymbol{e}}

\newcommand{\bq}{\boldsymbol{q}}

\newcommand{\bu}{\boldsymbol{u}}
\newcommand{\bv}{\boldsymbol{v}}
\newcommand{\bw}{\boldsymbol{w}}
\newcommand{\bx}{\boldsymbol{x}}

\newcommand{\bone}{\boldsymbol{1}}
\newcommand{\blambda}{\boldsymbol{\lambda}}

\newcommand{\bmu}{\boldsymbol{\mu}}

\newcommand{\BA}{\boldsymbol{A}}
\newcommand{\BB}{\boldsymbol{B}}
\newcommand{\BC}{\boldsymbol{C}}

\newcommand{\BH}{\boldsymbol{H}}

\newcommand{\BJ}{\boldsymbol{J}}
\newcommand{\BK}{\boldsymbol{K}}

\newcommand{\BM}{\boldsymbol{M}}
\newcommand{\BP}{\boldsymbol{P}}
\newcommand{\BQ}{\boldsymbol{Q}}
\newcommand{\BR}{\boldsymbol{R}}
\newcommand{\BS}{\boldsymbol{S}}
\newcommand{\BT}{\boldsymbol{T}}
\newcommand{\BU}{\boldsymbol{U}}
\newcommand{\BV}{\boldsymbol{V}}
\newcommand{\BW}{\boldsymbol{W}}
\newcommand{\BX}{\boldsymbol{X}}
\newcommand{\BY}{\boldsymbol{Y}}
\newcommand{\BZ}{\boldsymbol{Z}}

\newcommand{\BPhi}{\boldsymbol{\Phi}}

\newcommand{\BDelta}{\boldsymbol{\Delta}}
\newcommand{\BPi}{\boldsymbol{\Pi}}

\newcommand{\BLambda}{\boldsymbol{\Lambda}}
\newcommand{\BSigma}{\boldsymbol{\Sigma}}

\newcommand{\I}{\boldsymbol{I}}
\newcommand{\RR}{\mathbb{R}}
\newcommand{\lag}{\langle}
\newcommand{\rag}{\rangle}

\newcommand{\eps}{\epsilon}

\newcommand*\diff{\mathop{}\!\mathrm{d}}
\DeclareMathOperator{\Tr}{Tr}

\DeclareMathOperator{\VEC}{vec}

\DeclareMathOperator{\diag}{diag}
\DeclareMathOperator{\ddiag}{ddiag}

\DeclareMathOperator{\Corr}{corr}
\DeclareMathOperator{\Ran}{Ran}

\DeclareMathOperator{\mat}{mat}

\DeclareMathOperator{\SDR}{SDR}

\renewcommand{\Pr}{\mathbb{P}}

\long\def\\#1//{}

\definecolor{xl}{RGB}{200,50,120}

\begin{document}
\title{\bf On the Exactness of SDP Relaxation for \\ Quadratic Assignment Problem}
\author{Shuyang Ling\thanks{Shanghai Frontiers Science Center of Artificial Intelligence and Deep Learning, New York University Shanghai, Shanghai, China. S.L. and Z.S.Z. are (partially) financially supported by the National Key R\&D Program of China, Project Number 2021YFA1002800, National Natural Science Foundation of China (NSFC) No.12001372, Shanghai Municipal Education Commission (SMEC) via Grant 0920000112, and NYU Shanghai Boost Fund.}}

\maketitle

\begin{abstract}

Quadratic assignment problem (QAP) is a fundamental problem in combinatorial optimization and finds numerous applications in operation research, computer vision, and pattern recognition. However, it is a very well-known NP-hard problem to find the global minimizer to the QAP. In this work, we study the semidefinite relaxation (SDR) of the QAP and investigate when the SDR recovers the global minimizer. In particular, we consider the two input matrices satisfy a simple signal-plus-noise model, and show that when the noise is sufficiently smaller than the signal, then the SDR is exact, i.e., it recovers the global minimizer to the QAP. It is worth noting that this sufficient condition is purely algebraic and does not depend on any statistical assumption of the input data. We apply our bound to several statistical models such as correlated Gaussian Wigner model. Despite the sub-optimality in theory under those models, empirical studies show the remarkable performance of the SDR. Our work could be the first step towards a deeper understanding of the SDR exactness for the QAP.

\end{abstract}

\section{Introduction}

Given two matrices $\BA$ and $\BC$, how to find a simultaneous row and column permutation of $\BC$ such that the resulting two matrices are well aligned? This problem, known as the quadratic assignment (QAP), is one of most challenging problems in optimization~\cite{A00,A03,K72,L63,LDB+07}. Moreover, it has found numerous applications including graph matching~\cite{ABK15,CFS+04}, de-anonymization and privacy~\cite{OGE16}, and protein network~\cite{SXB08}, and traveling salesman~\cite{DS12,GW22}.

One of the most common approaches to find the optimal permutation is to minimize the least squares objective:
\begin{equation}\label{def:qap}
\min_{\BPi\in{\cal P}(n)}~\|\BA\BPi - \BPi\BC\|_F^2 ~~\Longleftrightarrow ~~\max~_{\BPi\in {\cal P}(n)}\lag \BA\BPi, \BPi\BC\rag
\end{equation}
where ${\cal P}(n)$ is the set of all $n\times n$ permutation matrices. 
In general, it is well-known NP-hard problem to find its global minimizer~\cite{K72}.
Numerous works have been done to either approximate or exactly solve the QAP problem~\cite{A00,A03,G62,L63,LDB+07}. Among various algorithms, convex relaxation is a popular approach to solve the QAP~\cite{DW09,FJBD13,LFF+16,PR09,ZKR+98,GW22}. In~\cite{G62}, Gilmore proposed a famous relaxation of quadratic assignment via linear programming and~\cite{G62,L63} derived a lower bound for the QAP problem.
Another straightforward convex relaxation is to relax the permutation matrix to the set of doubly stochastic matrix ${\cal D}(n)$, which leads to a quadratic program:
\begin{equation}\label{def:qpr}
\min_{\BX\in{\cal D}(n)}~\|\BA\BX - \BX\BC\|_F^2.
\end{equation}
One can also consider spectral relaxation of the QAP~\cite{RW92,A00,AB01}, and moreover an estimation of the optimal value of the QAP can be also characterized by the spectra of input data matrices.

In this work, we will focus on the semidefinite relaxation (SQR) of QAP. One well-known SDR was proposed in~\cite{ZKR+98}. After that, many variants of SDR have been proposed~\cite{PR09,DS10,DS12} to improve the formulation in~\cite{ZKR+98}, and~\cite{BKS18,OWX18} have studied efficient algorithms to solve the SDR. 
Our work will be more on the theoretical sides of the SDR for QAP. 
In particular, we are interested in the following question:
\begin{equation*}
\text{\em When does the SDR recover the global minimizer to~\eqref{def:qap}?}
\end{equation*}
Without any more constraints, there is certainly almost no hope to find the exact solution to the QAP for general input data due to its NP-hardness.
We consider a signal-plus-noise model for the quadratic assignment. More precisely, let $\BA$ be an $n\times n$ symmetric matrix (e.g. adjacency matrix), and $\BC$ is a perturbed matrix of $\BA$:
\begin{equation}\label{def:model}
\BC = \BPi^{\top}(\BA + \BDelta)\BPi
\end{equation}
where $\BPi\in {\cal P}(n)$ is an unknown permutation matrix. 
Our goal is to recover $\BPi$ from the two matrices $\BA$ and $\BC$ in an efficient way. 

In the noise-free case, i.e., $\BDelta= 0$ which corresponds to the graph isomorphism, the ground true permutation can be exactly recovered by solving~\eqref{def:qpr} or using spectral methods~\cite{FS15,KS18} under some regularity conditions on the spectra of $\BA$. Several other convex programs have been studied~\cite{ABK15}  to recover the ground true permutation for $\BDelta=0$ or in presence of very weak noise.
Recently, the quadratic assignment has been extensively studied under various statistical models (average-case analysis), especially in the context of graph matching or graph alignment such as
correlated Erd\"os-R\'enyi graph model~\cite{MRT23,WXY22,FMW+23b,HM23} and correlated Wigner model~\cite{CK16,G22,FMW+23}. 
A series of works fully exploit the statistical properties of random weight matrix, and use spectral methods~\cite{FQR19,FMW+23,GLM22} or extracts the feature of vertices~\cite{DMWX21,MRT21} to efficiently align two correlated random matrices. The core question is: under what noise levels, one can design an efficient algorithm to find the permutation~\cite{FMW+23,FMW+23b,FQR19} and whether the algorithm can achieve the information-theoretical threshold~\cite{CK16,G22}. 

On the other hand, the study on the optimization approaches to solve these random instances is quite limited compared with spectral methods or feature-extraction-based approaches. But optimization-based  approach often enjoys more robustness~\cite{CG18}. Therefore, we are interested in studying the performance of optimization methods, especially SDR, in solving the random instances of QAP. In the work~\cite{LFF+16}, the authors studied~\eqref{def:qpr} for the correlated  Erd\"os-R\'enyi model and proved that~\eqref{def:qpr} will never produce the exact permutation matrix even if the noise level is extremely small. The work~\cite{FMW+23,FMW+23b} proposed a spectral method to estimate the true permutation under both correlated Gaussian Wigner and Erd\"os-R\'enyi models. The spectral method can be viewed as a convex relaxation of~\eqref{def:qap}. Although that the global minimizer is not a permutation, it produces the true permutation after a rounding procedure.

SDR (semidefinite relaxation) has proven itself to be powerful in tacking many challenging nonconvex problems in signal processing and data science~\cite{VB96}.
The exactness or tightness of SDR has been studied in various problems in signal processing and data science. k-means and data clustering~\cite{ABCKVW15,TMPV17,LLL+20,LS20}, community detection~\cite{A17,ABH15}, synchronization~\cite{B18,ZB18,L22}, phase retrieval~\cite{CSV13}, matrix completion~\cite{R11}, and blind deconvolution~\cite{ARR13}. These works show that the SDR can recover the ground truth as long as the SNR (signal-to-noise ratio) is sufficiently large, such as the sample size is large enough or the noise in the data is smaller than certain threshold. Inspired by these observations, we will study when the exactness holds for the SDR of QAP, i.e., the SDR produces a permutation matrix which is also the global minimizer to~\eqref{def:qap}.

In this work, we focus on two variants of SDR for the QAP, and study its exactness under the signal-plus-noise model. We provide a sufficient condition that is based on the spectral gap of $\BA$ and also the noise strength $\|\BDelta\|$ to guarantee the exactness of the SDR. It is worth noting that this sufficient condition is deterministic, and can be applied to several statistical models. Despite the theoretical sub-optimality in the statistical examples, the SDR has shown powerful numerical performance and this could be the first step towards understanding the exactness of SDR for QAP.

\subsection{Notation}

Before proceeding, we go over some notations that will be used. We denote boldface $\bx$ and $\BX$ as a vector and matrix respectively, and $\bx^{\top}$ and $\BX^{\top}$ are their corresponding transpose. The $\I_n$ and $\BJ_n$ stands the $n\times n$ identity  and constant ``1" matrix of size $n\times n$. For a vector $\bx$, $\diag(\bx)$ is a diagonal matrix whose diagonal entries are given by $\bx.$
For two matrices $\BX$ and $\BY$ of the same size, $\BX\circ \BY$ is the Hadamard product, $\lag \BX,\BY\rag = \sum_{i,j}X_{ij}Y_{ij}= \Tr(\BX\BY^{\top})$ is their inner product, and $\BX\otimes\BY $ is their Kronecker product. For any matrix $\BX$, the Frobenius, operator norms, and the maximum absolute value among all entries are denoted by $\|\BX\|_F$, $\|\BX\|$, and $\|\BX\|_{\max}$ respectively. Here ${\cal P}(n)$ is  the set of $n\times n$ permutation matrix, ${\cal D}(n)$ is the set of all $n\times n$ doubly stochastic matrix, $\be_i$ is a one-hot vector, and $\delta_{ij} = 1$ if $i=j$, and $\delta_{ij} = 0$ if $i\neq j.$ We let $\VEC(\BX)$ be the vectorization of $\BX$ by stacking the columns of $\BX$, and for a given vector $\bx\in\RR^{n^2}$, we denote $\mat(\bx)\in\RR^{n\times n}$ as the matrization of $\bx$. 

\subsection{Organization}

The derivation of SDRs for the QAP is discussed in Section~\ref{s:main}; and we will present the main results and numerics in Section~\ref{s:main} as well. The theoretical justification of the main theorems will be given in Section~\ref{s:proof} and~\ref{s:proof2}.

\section{Preliminaries and main results}\label{s:main}

Without loss of generality, we assume the hidden ground true permutation is $\I_n$ in~\eqref{def:model}.
Note that we can rewrite this least squares objective by vectorizing $\BA\BPi$ and $\BPi\BC$:
\begin{align*}
\|\BA\BPi - \BPi\BC\|_F^2 & = \| (\I_n\otimes\BA - \BC\otimes\I_n)\bx\|^2 \\
& = \bx^{\top}  (\I_n\otimes\BA - \BC\otimes\I_n)^2\bx
\end{align*}
where $\BA$ and $\BC$ are symmetric, and $\bx$ is the vectorization of $\BPi$, i.e., for any permutation $\sigma(\cdot)$, which is bijective on $\{1,\cdots,n\}$, 
\begin{equation}\label{def:x}
\bx = 
\begin{bmatrix}
\be_{\sigma(1)} \\
\vdots \\
\be_{\sigma(n)}
\end{bmatrix}
\end{equation}
Therefore, by letting 
\begin{equation}\label{def:Mv1}
\BM: = (\I_n\otimes\BA - \BC\otimes\I_n)^2\succeq 0,
\end{equation}
the quadratic assignment problem is equivalent to
\begin{equation}\label{def:qaplift}
\min~\Tr(\BM\bx\bx^{\top}) \quad\text{ s.t. } \quad \bx = \VEC(\BPi),~~\BPi\in {\cal P}(n).
\end{equation}
Another equivalent form of~\eqref{def:qap} follows from 
\begin{align*}
\|\BA\BPi - \BPi\BC\|_F^2 & = \|\BA\|_F^2 - 2\lag \BA\BPi, \BPi\BC\rag + \|\BC\|_F^2  = -2 \lag (\I_n\otimes \BA)\bx, (\BC\otimes \I_n)\bx\rag + \|\BA\|_F^2 + \|\BC\|_F^2 \\
& = -2 \lag \BC\otimes \BA, \bx\bx^{\top}\rag + \|\BA\|_F^2 + \|\BC\|_F^2
\end{align*}
By letting 
\begin{equation}\label{def:Mv2}
\BM: = -\BC\otimes \BA,
\end{equation}
the least squares objective is equal to~\eqref{def:qaplift}.
Throughout the discussion, we will mainly focus on~\eqref{def:qaplift} with $\BM = (\I_n\otimes\BA - \BC\otimes \I_n)^2$ in~\eqref{def:Mv1}, and all the theoretical analysis also applies to~\eqref{def:Mv2} with minor modifications.

\subsection{Convex relaxation}

By letting $\BX = \bx\bx^{\top},$ we note that~\eqref{def:qaplift} is a linear function in $\BX$. Therefore,  the idea of convex relaxation of~\eqref{def:qaplift} is to find a proper convex set that includes all rank-1 matrix $\BX = \bx\bx^{\top}$ where $\bx = \VEC(\BPi)$ and $\BPi\in {\cal P}(n).$ By~\eqref{def:x},  $\BX = \bx\bx^{\top}$ is highly structured:
\[
\BX = 
\begin{bmatrix}
\be_{\sigma(1)}\be_{\sigma(1)}^{\top} & \be_{\sigma(1)} \be_{\sigma(2)}^{\top} & \cdots & \be_{\sigma(1)}\be_{\sigma(n)}^{\top}  \\
 \be_{\sigma(2)} \be_{\sigma(1)}^{\top} & \be_{\sigma(2)}\be_{\sigma(2)}^{\top} & \cdots &  \be_{\sigma(2)} \be_{\sigma(n)}^{\top} \\
 \vdots & \vdots & \ddots & \vdots \\
\be_{\sigma(n)}\be_{\sigma(1)}^{\top} & \be_{\sigma(n)} \be_{\sigma(2)}^{\top} & \cdots & \be_{\sigma(n)}\be_{\sigma(n)}^{\top}  \\
\end{bmatrix}
\]
which consists of $n^2$ blocks with each block exactly rank-1 and only containing one non-zero entry. 

Now we try to find a proper convex set which contains 
\[
\{\BX: \BX = \bx\bx^{\top}, \quad \bx = \VEC(\BPi),\quad \BPi\in{\cal P}(n)\}.
\]
It is obvious that $\BX\succeq 0$ and $\BX\geq 0$, which will be incorporated into the constraints.

\paragraph{Convex relaxation I.}

Note that for each $n\times n$ block, we have
\[
\BX_{ij} = \be_{\sigma(i)}\be_{\sigma(j)}^{\top}
\]
is exactly rank-1 and only contains one nonzero entry. 
Moreover, for any permutation $\sigma$, it holds that
\[
\Tr(\BX_{ij}) = 
\begin{cases}
1, & i =j, \\
0, & i\neq j,
\end{cases},~~~
\lag \BX_{ij}, \BJ_n\rag = 1
\]
and also
\begin{align*}
\sum_{i=1}^n \BX_{ii} & = \sum_{i=1}^n \be_{\sigma(i)}\be_{\sigma(i)}^{\top} = \I_n, \\
\sum_{i,j} \BX_{ij} & = \sum_{i,j}\be_{\sigma(i)}\be_{\sigma(j)}^{\top} = \sum_{i=1}^n \be_{\sigma(i)} \left(\sum_{j=1}^n \be_{\sigma(j)} \right)^{\top} = \BJ_n.
\end{align*}

Therefore, combining these constraints leads to the following convex relaxation:
\begin{equation}\label{def:sdr}
\begin{aligned}
\min\qquad& \lag \BM, \BX\rag \\
\text{s.t.} \qquad
& 
\BX\succeq 0, ~~~\BX \geq 0,  \\
& \Tr(\BX_{ij}) = \delta_{ij}, ~~~\lag \BX_{ij}, \BJ_n\rag = 1,~~\forall 1\leq i,j\leq n,  \\
& \sum_{i=1}^n \BX_{ii} = \I_n,~~~\sum_{i,j}\BX_{ij} = \BJ_n. 
\end{aligned}
\end{equation}

\paragraph{Convex relaxation II.} 

Note that for $\BX = \bx\bx^{\top}$, it holds
$\diag(\BX) = \bx.$ Using the fact that
\[
\begin{bmatrix}
\bx\bx^{\top} & \bx \\
\bx^{\top} & 1
\end{bmatrix} \succeq 0,
\]
we have a few new constraints:
\[
\begin{bmatrix}
\BX & \diag(\BX) \\
\diag(\BX)^{\top} & 1
\end{bmatrix} \succeq 0,~~~\mat(\diag(\BX))\in {\cal D}(n)
\]
where ${\cal D}(n)$ is the set of $n\times n$ doubly stochastic matrix.
Combining the constraints above with $\Tr(\BX_{ij}) = \delta_{ij}$ and nonnegativity, we get the SDR similar to~\cite{ZKR+98}:
\begin{equation}\label{def:sdr2}
\begin{aligned}
\min\qquad& \lag \BM, \BX\rag  \\
\text{s.t.} \qquad
& 
\begin{bmatrix}
\BX & \diag(\BX) \\
\diag(\BX)^{\top} & 1 
\end{bmatrix}\succeq 0, ~~~\BX \geq 0, \\
& \lag \BX_{ij}, \BJ_{n}\rag = 1,  ~~~ \Tr(\BX_{ij}) = \delta_{ij}, ~~\forall 1\leq i,j\leq n, \\
& \sum_{i=1}^n \BX_{ii} = \I_n, ~~~ \sum_{i,j}\BX_{ij} = \BJ_n,~~\forall 1\leq i,j\leq n,\\
& \mat(\diag(\BX))\bone_n = \mat(\diag(\BX))^{\top}\bone_n = \bone_n. 
\end{aligned}
\end{equation}
For the last constraints in~\eqref{def:sdr2}, the explicit form of $\mat(\diag(\BX))$ is given by
\[
\mat(\diag(\BX)) = 
\begin{bmatrix}
\BX_{11,11} & \BX_{11,22} & \BX_{11,33} & \cdots & \BX_{11,nn} \\
\BX_{22,11} & \BX_{22,22} & \BX_{22,33} & \cdots & \BX_{22,nn} \\
\BX_{33,11} & \BX_{33,22} & \BX_{33,33} & \cdots & \BX_{33,nn} \\
\vdots & \vdots & \vdots & \ddots & \vdots \\
\BX_{nn,11} & \BX_{nn,22} & \BX_{nn,33} & \cdots & \BX_{nn,nn}
\end{bmatrix}
\]
which reshapes the diagonal elements of $\BX$ into an $n\times n$ matrix. The relaxation~\eqref{def:sdr2} is tighter than~\eqref{def:sdr} as it imposes a few more constraints on  $\BX.$

\subsection{Main theorems}

With the introduction of two SDRs in~\eqref{def:sdr} and~\eqref{def:sdr2}, we are ready to present our main theorems.  Note that if $\BA$ has distinct eigenvalues and all eigenvectors are not orthogonal to $\bone_n$, then finding the optimal permutation is possible via simple convex relaxations~\cite{ABK15,KS18}. We interested in whether the exactness holds in presence of the noise, and below is our main theorem.
\begin{theorem}\label{thm:main}
Let $(\lambda_i,\bu_i)$ be the $i$-th eigenvalue and eigenvector of $\BA$, then $\BX= \bx\bx^{\top}$ is the unique global minimizer to SDR~\eqref{def:sdr} if 
\[
\min_{i\neq j}(\lambda_i - \lambda_j)^2  \min |\lag \bu_i,\bone_n\rag|^2 \geq n(\|\BDelta\|\|\BA\| + \|\BA\BDelta\|_{\max})
\]
where $\BA_i$ and $\BDelta_i$ are the $i$-th columns of $\BA$ and $\BDelta$ respectively.
\end{theorem}
There are remarks regarding Theorem~\ref{thm:main}. Suppose $\BDelta=0$, $\min_{i\neq j} |\lambda_i - \lambda_j| > 0$ and $\lag \bu_i,\bone_n\rag \neq 0$, then the exactness always holds. This is aligned with the results for graph isomorphism in~\cite{ABK15,KS18}. The interesting point is that even if the noise exists, as long as the noise strength is sufficiently small compared with minimum spectral gap and the alignment of $\bu_i$ and $\bone_n$, then SDR is still exact. In other words, our theorem provides a deterministic condition that guarantees the exactness of the SDR~\eqref{def:sdr}  in presence of noise. A result of similar flavor was derived in~\cite[Lemma 2]{ABK15}, which provides an error estimation between the solution from a quadratic program relaxation and the true permutation in terms of spectral gap and $\min_i |\lag \bu_i,\bone_n\rag|^2$. However, the exactness was not obtained in~\cite{ABK15}.

For the SDR in~\eqref{def:sdr2}, we can derive a similar deterministic sufficient condition for its exactness under the same assumptions.  
\begin{theorem}\label{thm:main2}
Let $(\lambda_i,\bu_i)$ be the $i$-th eigenvalue and eigenvector of $\BA$, then $\BX = \bx\bx^{\top}$ is the unique global minimizer to SDR~\eqref{def:sdr2}  if 
\[
\min_{i\neq j}(\lambda_i - \lambda_j)^2  \min |\lag \bu_i,\bone_n\rag|^2 \geq n(\|\BDelta\|\|\BA\| + \|\BA\BDelta\|_{\max})
\]
where $\BA_i$ and $\BDelta_i$ are the $i$-th columns of $\BA$ and $\BDelta$ respectively.
\end{theorem}

The proof of both theorems can be reduced to constructing a dual certificate (finding the dual variables) such that it can ensure the global optimality of the solution $\BX=\bx\bx^{\top}.$ This routine is well-established but the actual construction is highly problem dependent. For the SDR of QAP, the construction is not simple as the SDRs have complicated constraints.
One may expect a sharper theoretical bound for the SDR~\eqref{def:sdr2} than~\eqref{def:sdr}, as~\eqref{def:sdr2} has more constraints. However, due to the complication of constructing a proper dual certificate that ensures $\BX=\bx\bx^{\top}$, we employ a similar construction of the dual variables and this leads to the same deterministic condition.

As both Theorem~\ref{thm:main} and~\ref{thm:main2} are general, we present three special examples and see how our theorem works, and also make a comparison with the numerical experiments and the best known theoretical bounds. The experiments under these models imply that our theoretical bound is sub-optimal, although the SDR is remarkable in numerics. As briefly mentioned before, this sub-optimality results from  the construction of dual certificate which is quite challenging for the SDR of QAP with general input data. 

\subsection{Examples and numerics}
\paragraph{Example: diagonal matrix plus Gaussian noise.}
Suppose $\BA$ is a diagonal matrix, i.e., 
\[
\BA = \diag(\lambda_1,\cdots,\lambda_n)
\] 
where the eigenvalues are in a descending order and the eigenvector is $\bu_i = \be_i$. Suppose $\BDelta= \sigma\diag(\bw)$
where $\bw$ is a standard Gaussian random vector, and $\BC=\BA+\sigma\diag(\bw)$ is always a diagonal matrix. Then applying Theorem~\ref{thm:main} and~\ref{thm:main2} implies that the SDR is exact if
\[
\min_{i\neq j} (\lambda_i - \lambda_j)^2  \geq n\left( 4\sigma_{\SDR}\sqrt{\log n}\cdot \max_{1\leq i\leq n} |\lambda_i|\right)
\]
where $\|\BDelta\|\leq 2 \sqrt{\log n}$ holds with high probability and $\lag \be_i,\bone_n\rag = 1$.
Then the noise level $\sigma_{\SDR}$ should satisfy
\[
\sigma_{\SDR} \leq \frac{\min_{i\neq j} (\lambda_i - \lambda_j)^2}{4n\sqrt{\log n}\max_i |\lambda_i|}.
\]
On the other hand, given two diagonal matrices $\BA$ and $\BC$, the global minimizer to~\eqref{def:qap} should still be $\I_n$ if and only if the ordering of the eigenvalues remains unchanged, which holds if
\[
\sigma \leq \frac{\min_{i\neq j}|\lambda_i - \lambda_j|}{2\sqrt{\log n}}.
\]
For a specific example $\lambda_k = k$,~$1\leq k\leq n$, then
\[
\sigma_{\SDR} \leq \frac{1}{4n^2\sqrt{\log n}},~~~\sigma \leq \frac{1}{2\sqrt{\log n}}.
\]
This indicates that the SDR is sub-optimal by a factor of $n^2.$

We also look into the performance of~\eqref{def:sdr} in numerical simulations. We choose $n = 10$ due to the high computational complexity for larger $n$, and let $\lambda_k = k$ and $0\leq \sigma\leq 2$. For each $\sigma$, we run 20 experiments and compute the correlation of $\widehat{\BX}$ and $\bx\bx^{\top}:$
\begin{equation}\label{def:corr}
\Corr(\widehat{\BX}, \bx\bx^{\top}) = \frac{\bx^{\top}\widehat{\BX}\bx}{n^2}
\end{equation}
where $\widehat{\BX}$ is the solution to~\eqref{def:sdr}. This correlation is between 0 and 1: the higher the correlation, the better the recovery. In particular, if $\Corr(\cdot,\cdot)$ is 1, then $\widehat{\BX}= \bx\bx^{\top}$. We count one instance as exactness if $\Corr(\widehat{\BX}, \bx\bx^{\top}) \geq 1-10^{-3}.$
All simulations are done by using CVX~\cite{cvx} and the numerical result is presented below in Figure~\ref{fig:diag}.

\begin{figure}[h!]
\centering
\begin{minipage}{0.48\textwidth}
\centering
\includegraphics[width=75mm]{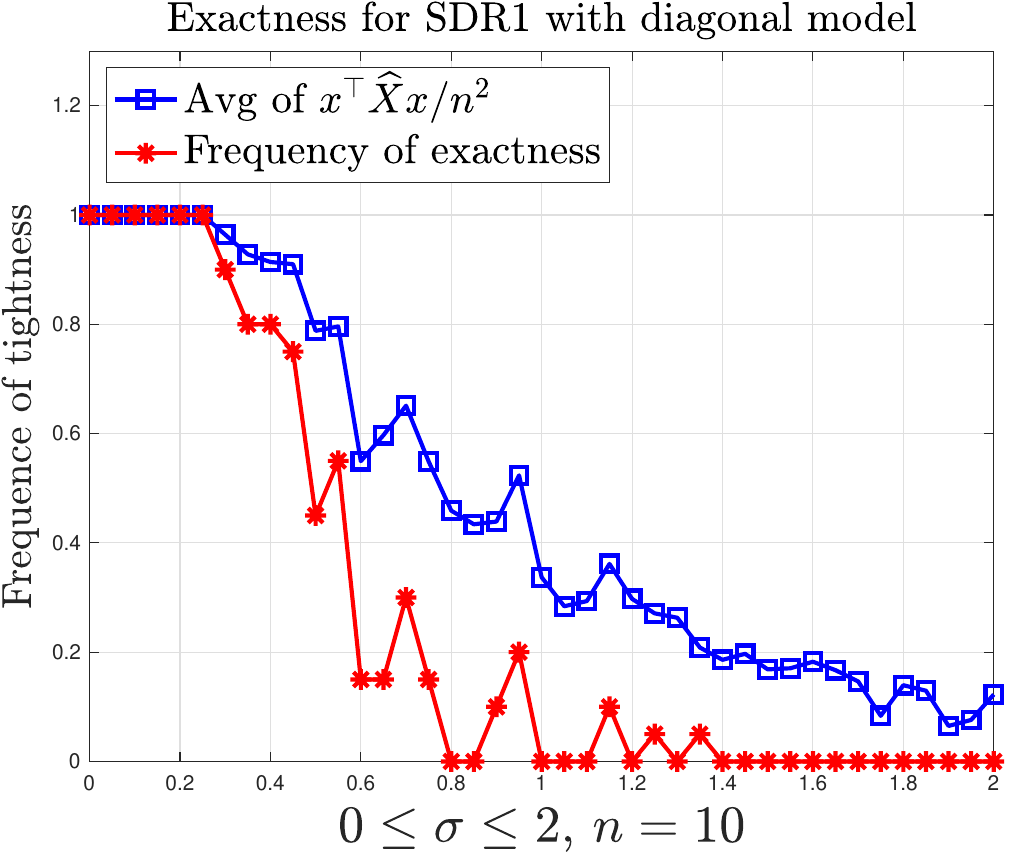}
\end{minipage}
\hfill
\begin{minipage}{0.48\textwidth}
\centering
\includegraphics[width=75mm]{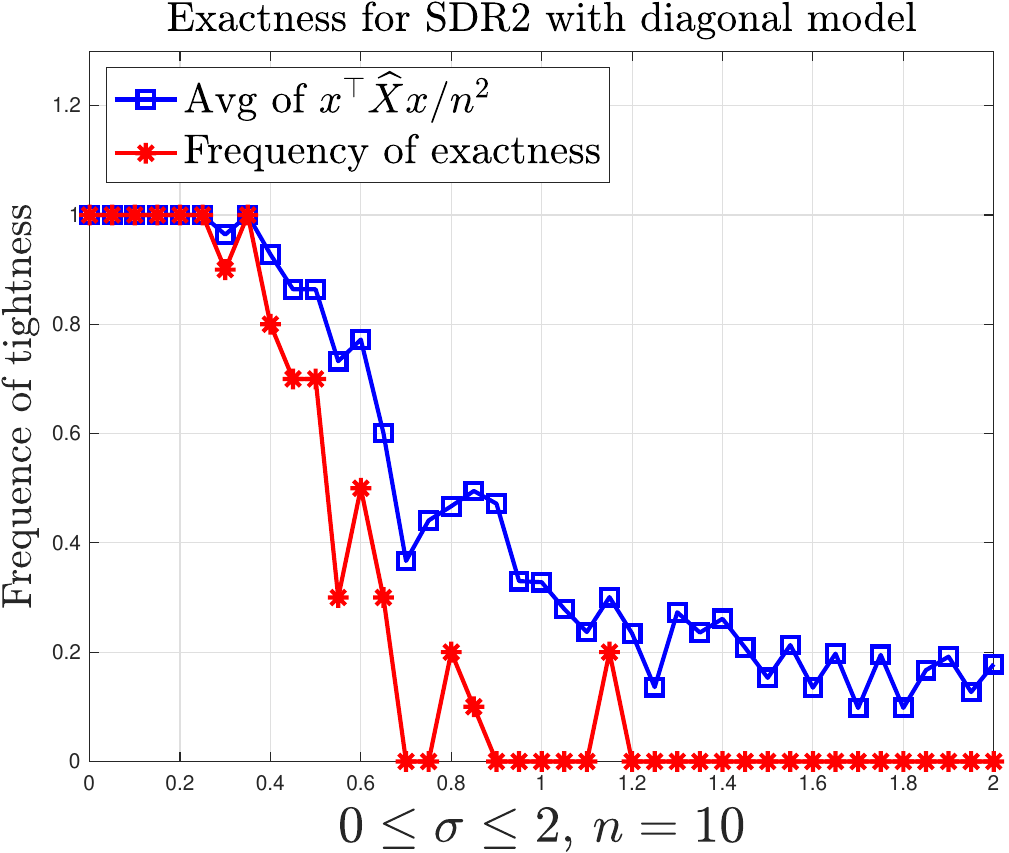}
\end{minipage}
\caption{Exactness of SDR~\eqref{def:sdr} and~\eqref{def:sdr2} for QAP under diagonal matrix model.}
\label{fig:diag}
\end{figure}
Figure~\ref{fig:diag} indicates that for $\sigma\leq 0.3$, then the SDR is exact and its performance nearly matches the optimal bound. The two types of SDRs perform similarly under this random setting.

\paragraph{Example: diagonal-plus-Wigner model}
Consider $\BA = \diag(\lambda_1,\cdots,\lambda_n)$ and $\BW$ is a Gaussian random matrix with $\diag(\BW) = 0$, and $\BC = \BA + \sigma\BW$.
Then applying Theorem~\ref{thm:main} implies that
\[
\min_{i\neq j} (\lambda_i - \lambda_j)^2 \min |\lag \bu_i,\bone_n\rag|^2 \geq n \sigma_{\SDR} \max_i |\lambda_i| \cdot 2\sqrt{n}~~\Longleftrightarrow~~
\sigma_{\SDR} \leq \frac{\min_{i\neq j} (\lambda_i - \lambda_j)^2 }{2n^{3/2}\max_i|\lambda_i|}
\]
is needed to ensures the exactness where $\bu_i = \be_i$ and $\|\BW\|\leq (2 +o(1))\sqrt{n}$ holds with high probability~\cite{V18}. In particular, for $\lambda_k = k$,~$1\leq k\leq n$, then the exactness of the SDR holds if $\sigma_{\SDR} \leq 1/(2n^{5/2}).$

\begin{figure}[h!]
\centering
\begin{minipage}{0.48\textwidth}
\centering
\includegraphics[width=75mm]{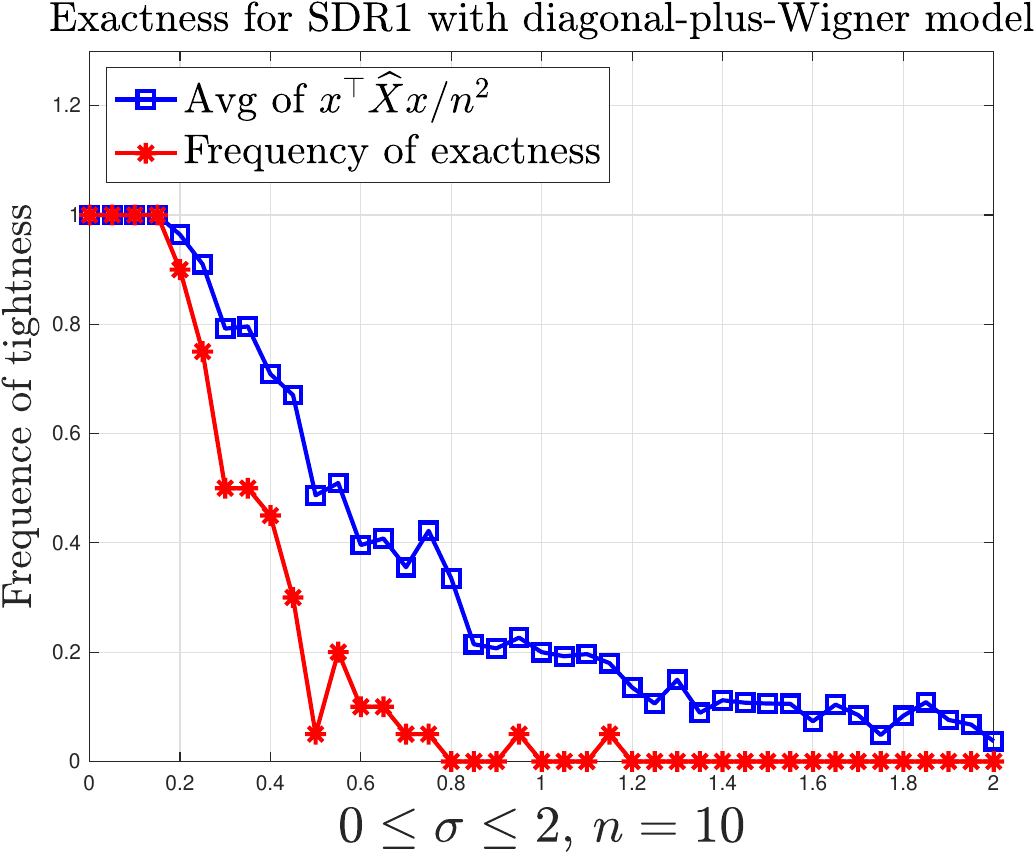}
\end{minipage}
\hfill
\begin{minipage}{0.48\textwidth}
\centering
\includegraphics[width=75mm]{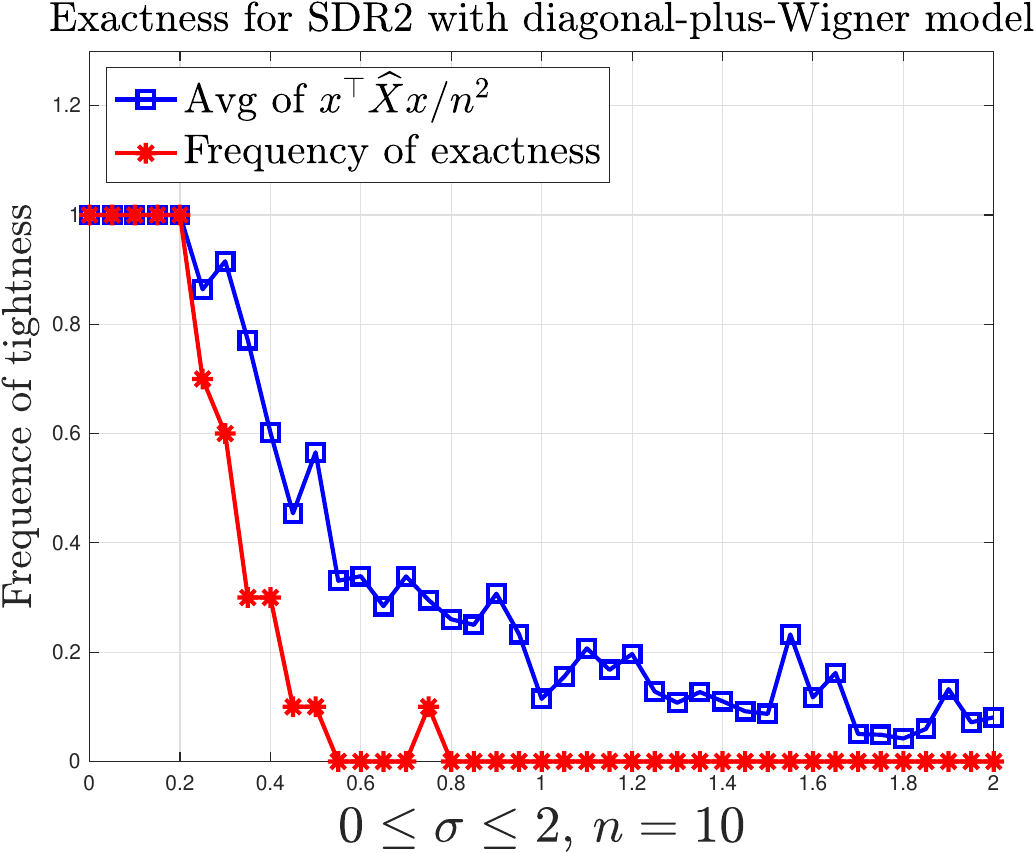}
\end{minipage}
\caption{Exactness of SDR~\eqref{def:sdr} and~\eqref{def:sdr2} for QAP under diagonal-plus-Wigner matrix model.}
\label{fig:diaggaussian}
\end{figure}
The numerical experiment is given in Figure~\ref{fig:diaggaussian}, which shows the exactness holds if $\sigma \leq 0.15$. For $n = 10$, our bound requires $\sigma_{\SDR}\leq 0.0016$ in theory. This implies that there is gap between the actual performance of the SDR and our results. 

\paragraph{Example: correlated Wigner model.} 
The correlated Wigner model assumes $\BA$ and $\BC$ are two Gaussian random matrices satisfying 
$\BC = \BA + \sigma\BW$
where $\BW$ is a Gaussian random matrix independent of $\BA$. To apply Theorem~\ref{thm:main} and~\ref{thm:main2}, it suffices to obtain a lower bound on the spectral gap and also an upper bound on $\|\BA\|.$ 

\begin{figure}[h!]
\centering
\begin{minipage}{0.48\textwidth}
\centering
\includegraphics[width=75mm]{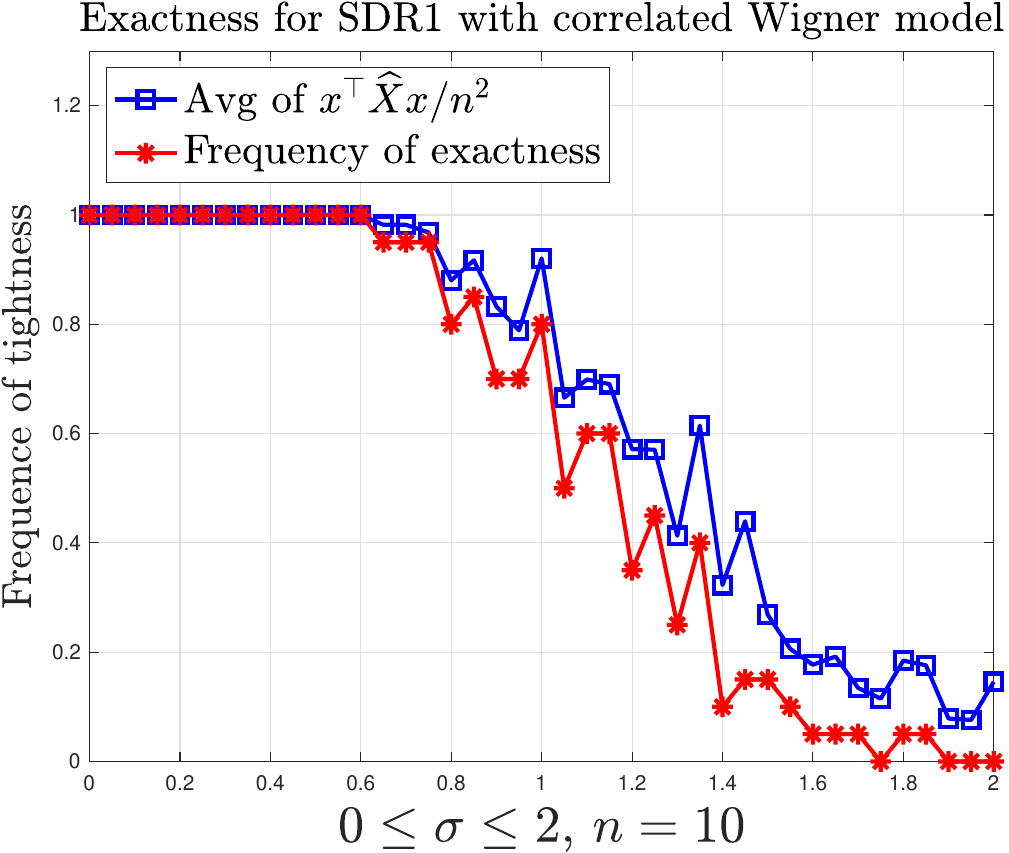}
\end{minipage}
\hfill
\begin{minipage}{0.48\textwidth}
\centering
\includegraphics[width=75mm]{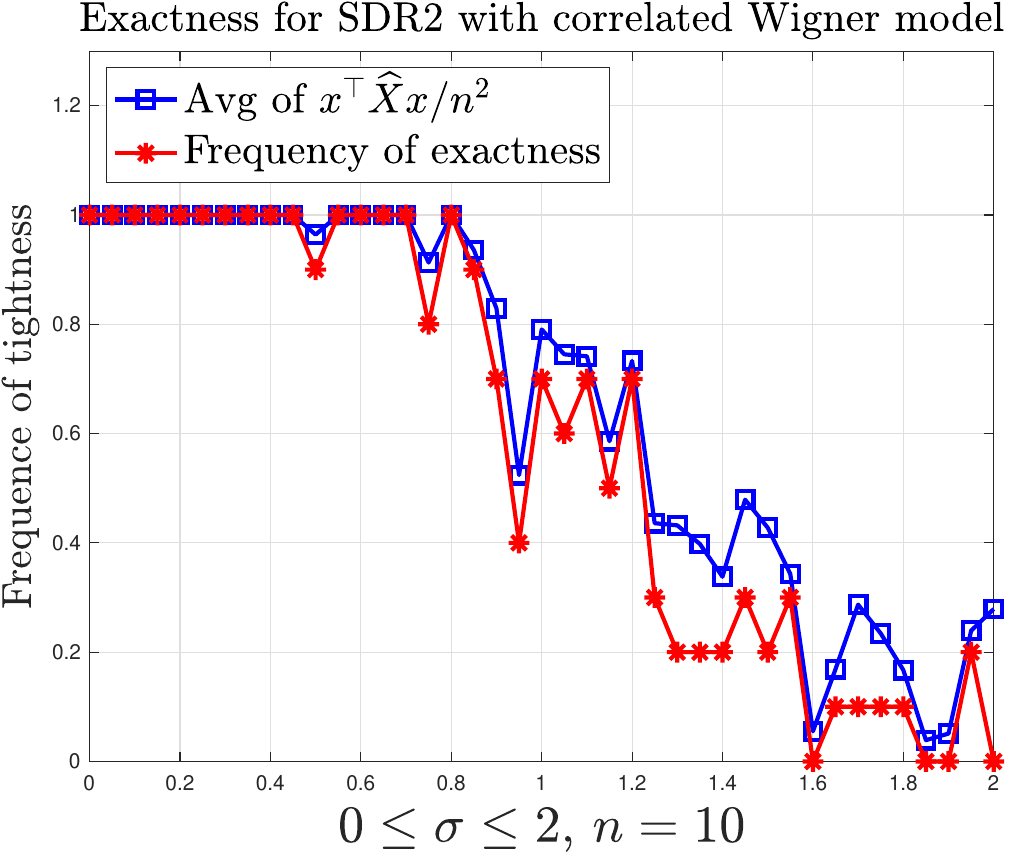}
\end{minipage}
\caption{Exactness of SDR~\eqref{def:sdr} and~\eqref{def:sdr2} for QAP under correlated-Wigner model.}
\label{fig:wigner}
\end{figure}

For the spectral gap,~\cite[Corollary 1]{FTW19} implies for $n$ sufficiently large, it holds that
\[
\Pr(2^{-3/2}n \delta_n \geq t) \approx 2\int_t^{\infty} x e^{-x^2}\diff x = e^{-t^2} \geq 1-t^2
\]
where $\delta_n$ denotes the smallest gap.
Therefore, it holds with probability at least $1-O(\log^{-1}n)$ that
\[
\delta_n = \max_{i\neq j}|\lambda_i - \lambda_j| \geq \frac{1}{n\sqrt{\log n}}.
\]
As a result, our result implies that the exactness holds if
\[
\frac{1}{2}\left(\frac{1}{n\sqrt{\log n}} \right)^2 \geq n (2\sigma\sqrt{n}\cdot 2\sqrt{n})
\Longrightarrow \sigma_{\SDR}\lesssim \frac{1}{n^4\log n}
\]
where $\|\BW\|\leq 2(1+o(1))\sqrt{n}$ and $\min_i|\lag \bu_i,\bone_n\rag|^2 \geq 1/2$ holds with high probability.
As shown in the Figure~\ref{fig:wigner}, it implies that $\sigma= 0.6$ suffices to have the exactness numerically which is far better than the theoretic bound. Under the correlated Wigner model,~\eqref{def:sdr2} performs better than~\eqref{def:sdr} when $\sigma\geq 0.6.$

Let's conclude this section by commenting on the numerics and also future research problems. In all three examples, we can see the SDRs perform much better in numerics than the theoretical bound by a factor of $n^2$ or even more. In fact, for the correlated Gaussian Wigner model, the SDRs are able to recover the true permutation even if $\sigma$ is of constant order, which is comparable with the state-of-the-art works on the same model such as~\cite{FMW+23}. Therefore, there is still big room for the improvement. In particular, it will be very interesting to show theoretically that SDR achieves exactness on random instances that is consistent with the numerical observation. This asks for a new technique to construct dual certificates that could exploit the statistical property of noise. More importantly, it is unclear that for the QAP problem, what is the best way to characterize the SNR (signal-to-noise ratio) that could be useful in analyzing the SDR of QAP. For example, in the correlated Gaussian Wigner model, using the spectral gap for the analysis does not seem to match the numerical performance as the spectral gap of a Gaussian random matrix can be very small for larger $n$. All of these will be worthwhile to consider in the future. 

\section{Proof of Theorem~\ref{thm:main}}\label{s:proof}

The proof of Theorem~\ref{thm:main} essentially follows from the well-established technique in convex relaxation. Without loss of generality, we assume $\BPi = \I_n$ and thus
\[
\bx^{\top} = [\be_1^{\top},\cdots,\be_n^{\top}]\in\RR^{n^2},~~\|\bx\|^2 = n.
\]
Our goal is to establish sufficient conditions that ensure the global minimizer to~\eqref{def:sdr} is $\bx\bx^{\top}$ by constructing a dual certificate.
 
\subsection{Dual program and optimality condition of~\eqref{def:sdr}}

For each constraint in~\eqref{def:sdr}, we assign a dual variable:
\begin{align*}
\BQ \succeq 0 &: \BX\succeq 0, \\
\BB \geq 0 &: \BX \geq 0, \\
\BT\in\RR^{n\times n} &: \Tr(\BX_{ij}) = \delta_{ij}, \\
\BZ\in\RR^{n\times n} &: \lag \BX_{ij},\BJ_n\rag = 1, \\
\BK \in\RR^{n\times n} &: \sum_{i=1}^n \BX_{ii} = \I_n, \\
\BH\in \RR^{n\times n} &: \sum_{i,j} \BX_{ij} = \BJ_n.
\end{align*} 
The Lagrangian function is
{\small
\begin{align*}
& {\cal L}(\BX,\BQ,\BB,\BT,\BZ,\BK,\BH) \\
&  \qquad  : = \lag \BX,\BM-\BQ -\BB  ) \rag  - \sum_{i,j} t_{ij}(\Tr(\BX_{ij}) -\delta_{ij})  - \sum_{i,j} z_{ij}(\lag \BX_{ij},\BJ_n\rag - 1)   \\
& \qquad\qquad -\left \lag \sum_{i=1}^n \BX_{ii} - \I_n, \BK\right\rag - \left\lag \sum_{i,j}\BX_{ij} - \BJ_n, \BH\right\rag \\
&  \qquad  = \lag  \BX,\BM-\BQ-\BB - \BT\otimes \I_n - \BZ\otimes \BJ_n -  \I_n\otimes \BK -  \BJ_n\otimes \BH\rag  + \Tr(\BT + \BK) + \lag \BJ_n,\BZ+\BH\rag.
\end{align*}
}

As a result, the dual program of~\eqref{def:sdr} is
\begin{equation}\label{eq:dual}
\begin{aligned}
\max\qquad &  \Tr(\BT + \BK) + \lag \BJ_n,\BH+\BZ\rag \\
\text{s.t.}\qquad
& \BM = \BQ+\BB+ \BT\otimes \I_n + \BZ\otimes \BJ_n +  \I_n\otimes \BK + \BJ_n\otimes \BH,  \\
& \BQ \succeq 0,~~  \BB\geq 0. 
\end{aligned}
\end{equation} 

Suppose $\BX = \bx\bx^{\top}$ is the global minimizer, then corresponding the complementary slackness: 
\begin{equation}\label{eq:compsl}
\lag \BB, \BX\rag = 0,\quad\BQ\BX = 0.
\end{equation}
Note that $\BQ\succeq 0$ and $\BX = \bx\bx^{\top}\succeq 0$, and then $\BQ\BX=0$ is equivalent to
$\BQ\bx = 0.$
Therefore the KKT condition becomes:
\begin{enumerate}
\item Stationarity:
\begin{equation}\label{eq:kkt1}
\BM = \BQ+\BB+ \BT\otimes \I_n +  \I_n\otimes \BK + \BJ_n\otimes \BH+ \BZ\otimes \BJ_n,
\end{equation}

\item Dual feasibility:
\begin{equation}\label{eq:kkt2}
\BB \geq 0,\qquad \BQ  \succeq 0,
\end{equation}

\item Complementary slackness:
\begin{equation}\label{eq:kkt3}
\lag \BB, \BX\rag = 0,\qquad \BQ\bx = 0,
\end{equation}
\end{enumerate}
where $\BM = (\I_n\otimes\BA - \BC\otimes\I_n)^2$ is the data matrix. Due to the nonnegativity of $\BB$, $\lag \BB,\BX\rag = 0$ is equivalent to $\BB\circ \bx\bx^{\top} = 0$ where ``$\circ$" is the Hadamard product of two matrices.

Now we compute $\BQ\bx = 0$ to further make the KKT condition explicit. Instead of computing $\BQ\bx$, we will use the property of Kronecker product to simplify the expression. Note that
\begin{align*}
&  \mat(\BM\bx) = \mat((\I_n\otimes\BA - \BC\otimes\I_n)^2\bx) \\ 
& \qquad\qquad~~ = \BA^2 - 2\BA\BC +\BC^2 = \BDelta^2 + \BDelta\BA - \BA\BDelta, \\
&  \mat((\BT\otimes\I_n)\bx) = \BT, \\
&  \mat((\BZ\otimes\BJ_n)\bx) = \BJ_n\BZ, \\
& \mat(( \I_n\otimes \BK)\bx) = \BK, \\
& \mat((\BJ_n\otimes\BH)\bx) = \BH\BJ_n,
\end{align*}
where $\bx = \VEC(\I_n).$
We define
\begin{equation}\label{def:S}
\BS =  \BT\otimes \I_n +  \I_n\otimes \BK + \BJ_n\otimes \BH+ \BZ\otimes \BJ_n
\end{equation}
as we will use it quite often, and it holds
\begin{align*}
\mat(\BS\bx) & = \BT + \BJ_n\BZ + \BK + \BH\BJ_n, \\
\mat((\BM-\BS)\bx) & = \BDelta^2 + \BDelta\BA - \BA\BDelta - (\BT + \BJ_n\BZ + \BK + \BH\BJ_n).
\end{align*}
Therefore,~\eqref{eq:kkt1},~\eqref{eq:kkt2} and~\eqref{eq:kkt3} imply
\begin{equation}\label{eq:MSx}
\begin{aligned}
\mat((\BQ+\BB)\bx)  & = \mat((\BM-\BS)\bx)\\
& = \mat( (\BM -( \BT\otimes \I_n +  \I_n\otimes \BK + \BJ_n\otimes \BH+ \BZ\otimes \BJ_n))\bx) \\
& = \BDelta^2 +\BDelta\BA-\BA\BDelta -( \BT+\BJ_n\BZ + \BK+\BH\BJ_n) = \mat(\BB\bx)
\end{aligned}
\end{equation}
where $\BQ\bx=0$ is in~\eqref{eq:kkt2}.
Therefore, the KKT condition has a simplified form:
\begin{equation}\label{eq:kktfull}
\begin{aligned}
& \BQ = \BM-\left(\BB+ \BT\otimes \I_n + \BZ\otimes \BJ_n+  \I_n\otimes \BK + \BJ_n\otimes \BH\right)\succeq 0, \\
&  \BDelta^2 + \BDelta\BA - \BA\BDelta =  \mat(\BB\bx)  + (\BT +\BJ_n\BZ+\BK +\BH\BJ_n ), \\
& \lag \BB, \BX\rag = 0,~~\BB \geq 0.
\end{aligned}
\end{equation}
In the following section, we will construct an explicit dual certificate, i.e., $\BT, \BZ,\BK,\BH$, and $\BB\succeq 0$ such that~\eqref{eq:kktfull} holds. Moreover, we can also certify $\BX=\bx\bx^{\top}$ as the unique global minimizer under mild conditions. 

\begin{theorem}\label{thm:suff}
Suppose there exist $\BB\geq 0$ with $\BB\circ\bx\bx^{\top} = 0$, and $\BT,\BK,\BZ$ and $\BH$ such that $\BQ\succeq c\left(\I_{n^2} - n^{-1}\bx\bx^{\top}\right)$ for some $c\geq 0$, then $\BX= \bx\bx^{\top}$ is a global minimizer. In addition, if $c>0$, then $\BX$ is the unique global minimizer.
\end{theorem}
\begin{proof}
Let $\widehat{\BX}$ be any feasible solution that is not equal to $\bx\bx^{\top}$. Then 
\begin{align*}
\lag \BQ,\widehat{\BX}\rag & = \lag \BQ, \widehat{\BX} - \bx\bx^{\top}\rag = \lag \BM - \BS-\BB,  \widehat{\BX} - \bx\bx^{\top}\rag \\
& = \lag \BM -\BB,  \widehat{\BX} - \bx\bx^{\top}\rag \\
& = \lag \BM, \widehat{\BX} - \bx\bx^{\top}\rag - \lag \BB, \widehat{\BX}\rag
\end{align*}
where $\BQ = \BM - \BS- \BB$ and $\lag \BS, \widehat{\BX} - \bx\bx^{\top} \rag = 0$ holds due to the primal feasibility of $\widehat{\BX}$ and $\bx\bx^{\top}$. Now we have
\begin{align*}
 \lag \BM, \widehat{\BX} - \bx\bx^{\top}\rag & = \lag \BQ+\BB,\widehat{\BX}\rag  \geq \lag \BQ, \widehat{\BX}\rag= \left\lag \BQ, \left(\I_{n^2} - \frac{\bx\bx^{\top}}{n}\right)\widehat{\BX}\left(\I_{n^2} - \frac{\bx\bx^{\top}}{n}\right)\right\rag \geq 0
\end{align*}
where $ \left(\I_{n^2} - n^{-1}\bx\bx^{\top}\right)\widehat{\BX}\left(\I_{n^2} - n^{-1}\bx\bx^{\top}\right)$ is nonzero. In particular, if $c> 0$, then the inequality above is strict, implying the uniqueness of the global minimizer. 
\end{proof}

\subsection{Construction of a dual certificate}

From the first equation in~\eqref{eq:kktfull} and $\BQ\bx = 0$, it holds
\[
\BB\bx = (\BM - \BS)\bx \geq 0
\]
where $\BB\geq 0$ and $\BS = \BT\otimes \I_n + \BZ\otimes \BJ_n+  \I_n\otimes \BK + \BJ_n\otimes \BH.$
Note that $\BB\circ \bx\bx^{\top} = 0$, and thus $\diag(\bx)\BB\bx = 0$, i.e., the diagonal elements of $\mat(\BB\bx)$ equal 0. Therefore, the second equation in~\eqref{eq:kktfull} gives
$\diag(\BDelta^2) = \diag(\BT+\BK) + (\BZ+\BH)\bone_n$.

With the discussion above, we try to first determine $\BB$ via solving the linear equation $\BB\bx = (\BM-\BS)\bx$, and we will find that the construction of a dual certificate finally reduce to searching for proper $\BT,\BK,\BZ$ and $\BH$.
\begin{proposition}[\bf Construction of $\BB$ and $\BQ$]\label{prop:consBQ}
Suppose 
\begin{equation}\label{eq:Bxdiag}
\diag(\BDelta^2) = \diag(\BT+\BK) + (\BZ+\BH)\bone_n
\end{equation}
for some $\BT,\BK$, $\BZ$, and $\BH$. Let
\begin{equation}\label{def:B}
\BB =  \frac{1}{n}\left(\BM -\BS\right)\bx\bx^{\top}  + \frac{1}{n}\bx\bx^{\top}\left(\BM -\BS\right),
\end{equation}
then $\BB\bx = (\BM-\BS)\bx$ and $\BB\circ \bx\bx^{\top} = 0$ hold automatically. Moreover, $\BB\geq 0$ if
\begin{equation}\label{def:MSx}
(\BM-\BS)\bx\geq 0 \Longleftrightarrow \BDelta^2 +\BDelta\BA-\BA\BDelta - ( \BT + \BJ_n\BZ + \BK + \BH\BJ_n) \geq 0.
\end{equation}
for some $\BT,\BZ,\BK$ and $\BH.$ As a result, it holds
\begin{equation}\label{def:Qcons}
\BQ = \BM - \BS-\BB= \left( \I_{n^2} - \frac{\bx\bx^{\top}}{n}\right) \left(\BM -\BS\right)\left( \I_{n^2} - \frac{\bx\bx^{\top}}{n}\right) 
\end{equation}
with the construction of $\BB$ in~\eqref{def:B}.
\end{proposition}
\begin{proof}
Consider $\BB$ in~\eqref{def:B}, it holds that
\begin{align*}
\BB \bx 
& = \frac{1}{n}\left(\BM -\BS\right)\bx\bx^{\top}\bx  + \frac{1}{n}\bx\cdot \bx^{\top}\left(\BM -\BS\right)\bx = \left(\BM -\BS\right)\bx 
\end{align*}
where 
\begin{align*}
\bx^{\top}(\BM-\BS)\bx& = \lag \I_n, \mat((\BM-\BS)\bx) \rag \\
& = \lag \I_n, ( \BDelta^2 +\BDelta\BA-\BA\BDelta -\mat(\BS\bx)) \rag \\
& = \Tr\left( \diag(\BDelta^2)- ( \diag(\BT+\BK) + (\BH+\BZ)\bone_n)\right) = 0
\end{align*}
and the expression of $(\BM-\BS)\bx$ follows from~\eqref{eq:MSx} and~\eqref{eq:Bxdiag}.

We proceed to verify that $\BB\circ\bx\bx^{\top} = 0$:
\begin{align*}
\mat(\diag(\bx)\BB\bx) &=  \mat(\diag(\bx)(\BM-\BS)\bx) \\
& = \I_n\circ \mat((\BM-\BS)\bx) \\
& = \I_n \circ ( \BDelta^2 +\BDelta\BA-\BA\BDelta - ( \BT + \BJ_n\BZ + \BK + \BH\BJ_n)) \\
& = \diag(\BDelta^2) - (\diag(\BT +\BK) + (\BZ+\BH)\bone_n) = 0
\end{align*}
follows from~\eqref{eq:Bxdiag}. This implies that the diagonal entries of $\mat(\BB\bx)$ are zero, i.e., the supports of $\BB\bx$ and $\bx$ are disjoint.
Now using $\BB$ in~\eqref{def:B} leads to 
\begin{align*}
\BB\circ \bx\bx^{\top}& = \diag(\bx)\BB\diag(\bx) = \frac{1}{n}\diag(\bx)\left(\BM -\BS\right)\bx\bx^{\top}  + \frac{1}{n}\bx\bx^{\top}\left(\BM -\BS \right)\diag(\bx) \\
& = \frac{1}{n}\diag(\bx)\BB \bx\bx^{\top}  + \frac{1}{n}\bx\bx^{\top}\BB\diag(\bx) =0
\end{align*}
where $\diag(\bx) \bx = \bx$ and $\diag(\bx)\BB\bx = 0.$


Finally, to ensure $\BB\geq 0$, it suffices to have
\[
(\BM-\BS)\bx = \BDelta^2 +\BDelta\BA-\BA\BDelta - ( \BT + \BJ_n\BZ + \BK + \BH\BJ_n) \geq 0
\]
due to the construction of $\BB$ in~\eqref{def:B}.

Then it holds that
\begin{align*}
\BQ & = \BM - \BS - \BB \\
& =\BM - \BS - \frac{1}{n}\left(\BM -\BS\right)\bx\bx^{\top}  - \frac{1}{n}\bx\bx^{\top}\left(\BM -\BS\right)\left( \I_{n^2} - \frac{\bx\bx^{\top}}{n}\right) \\
& = \left( \I_{n^2} - \frac{\bx\bx^{\top}}{n}\right) \left(\BM -\BS\right)\left( \I_{n^2} - \frac{\bx\bx^{\top}}{n}\right) 
\end{align*}
where $\bx^{\top}(\BM-\BS)\bx = 0.$
\end{proof}

Finally, we summarize our findings: by choosing $\BB$ in the form of~\eqref{def:B}, it suffices to find $\BT$, $\BK$, $\BZ$ and $\BH$ such that 
\begin{equation}\label{eq:kktfinal_simp}
\begin{aligned}
& \diag(\BT+\BK) + (\BZ+\BH)\bone_n = \diag(\BDelta^2), \\
&  \BDelta^2 +\BDelta\BA-\BA\BDelta - ( \BT + \BJ_n\BZ + \BK + \BH\BJ_n) \geq 0, \\
&  \left( \I_{n^2} - \frac{\bx\bx^{\top}}{n}\right) \left(\BM -( \BT\otimes \I_n + \BZ\otimes \BJ_n+  \I_n\otimes \BK + \BJ_n\otimes \BH)\right)\left( \I_{n^2} - \frac{\bx\bx^{\top}}{n}\right)  \succeq 0.
\end{aligned}
\end{equation}
where the first two constraints ensure $\BB\circ\bx\bx^{\top} = 0$ and $\BB\geq 0$ respectively, and third one corresponds to $\BQ\succeq 0.$

Now we proceed to prove Theorem~\ref{thm:main}. The proof relies on the following proposition that assume $\BDelta = 0$.
In the noiseless case, i.e., $\BC = \BA$ with $\BDelta = 0$:
\begin{align*}
\BM & = (\I_n\otimes\BA - \BA\otimes\I_n )^2 =  \sum_{i=1}^n\sum_{j=1}^n (\lambda_i  - \lambda_j)^2\bu_i\bu_i^{\top}\otimes\bu_j\bu_j^{\top}
\end{align*}
where $(\lambda_i,\bu_i)$ is the $i$th eigenvalue and eigenvector of $\BA.$ 
The smallest eigenvectors of $\BM$ are 0 with corresponding eigenvectors $\{\bu_i\otimes\bu_i,1\leq i\leq n\}.$
\begin{theorem}[\bf Noise-free version of Theorem~\ref{thm:main}]\label{thm:noisefree}
Suppose $|\lag \bu_i,\bone_n\rag|^2> 0$ for all $1\leq i\leq n$ and the eigenvalues of $\BA$ are distinct. Moreover, by letting
\begin{equation}\label{eq:choiceTKHZ}
\BT =-t\BJ_n,~\BK =-t\BJ_n, ~\BH = t\BJ_n/n,~\BZ = t\BJ_n/n,
\end{equation}
then the second smallest eigenvalue of 
\[
\BQ = \BM - \BS = \sum_{i=1}^n\sum_{j=1}^n (\lambda_i  - \lambda_j)^2\bu_i\bu_i^{\top}\otimes\bu_j\bu_j^{\top} +  t \BJ_n\otimes(\I_n - \BJ_n/n) + t(\I_n - \BJ_n/n)\otimes\BJ_n.
\] 
is bounded below by
\[
\lambda_2(\BQ) \geq \frac{2}{n} \min_{i\neq j}(\lambda_i - \lambda_j)^2 \cdot \min_{1\leq i\leq n} |\lag \bu_i,\bone_n\rag|^2
\]
for sufficiently large $t>0$, and $\BQ\bx = 0.$
\end{theorem}
\begin{proof}
Under~\eqref{eq:choiceTKHZ}, it holds
\begin{align*}
\BS & = -t \left( \BJ_n\otimes(\I_n - \BJ_n/n) + t(\I_n - \BJ_n/n)\otimes\BJ_n\right), \\
\BQ & = \BM - \BS = \BM +  t \BJ_n\otimes(\I_n - \BJ_n/n) + t(\I_n - \BJ_n/n)\otimes\BJ_n.
\end{align*}
It is easy to verify that
$\BM\bx = \BS\bx = 0$, and thus $\BQ\bx = 0$ and $\BB=0.$

Next, we will show $\lambda_2(\BM-\BS)> 0$ under the assumption of this theorem. 
Before proceeding, we introduce a few notations.
Let $\BU = [\bu_1,\cdots,\bu_n]$ consist of $n$ eigenvectors of $\BA$ with corresponding eigenvalues $\BLambda = \diag(\lambda_1,\cdots,\lambda_n)$, i.e., $\BA = \BU\BSigma\BU^{\top}.$
Let
\[
\BPhi = [\bu_1\otimes\bu_1,\cdots,\bu_n\otimes\bu_n]\in\RR^{n^2\times n}
\]
be the eigenvectors of $\BM$ w.r.t. eigenvalue $0$.
Then it holds that
\begin{equation}\label{eq:phiphix}
\BPhi\BPhi^{\top} \bx = \sum_{i=1}^n (\bu_i\bu_i^{\top}\otimes\bu_i\bu_i^{\top})\bx = \bx,~~~~\BPhi^{\top}\bx = \bone_n.
\end{equation}
Let $\BP = \BJ_n/n$, and 
\begin{align*}
\bar{\BP} & = (\BP\otimes\I_n - \I_n \otimes\BP)^2 = \frac{1}{n}\left( \BJ_n\otimes \left(\I_n - \frac{\BJ_n}{n}\right) + \left(\I_n - \frac{\BJ_n}{n}\right)\otimes \BJ_n\right)~~~\BPi = \I_{n^2} - \frac{\bx\bx^{\top}}{n},
\end{align*}
where $\bar{\BP}$, $\BP$, $\BPi,$ and $\BPi - \bar{\BP}$ are all projection matrices satisfying $\BPi \bar{\BP} = \bar{\BP}$ and $\bar{\BP} \bx = 0.$

Our goal is to show that for some sufficiently large $t>0$, it holds
\begin{align*}
\BM - \BS &= \BM+ t  \BJ_n\otimes(\I_n - \BJ_n/n) + t(\I_n - \BJ_n/n)\otimes\BJ_n = \BM+ nt\bar{\BP} \succeq c\BPi
\end{align*}
for some $c>0$. 
Lemma~\ref{lem:supp} implies that it suffices to prove that
\[
(\BPi - \bar{\BP}) \BM(\BPi - \bar{\BP})\succ c (\BPi - \bar{\BP}).
\]
Note that 
\[
\BM \succeq \min_{i\neq j}(\lambda_i - \lambda_j)^2 \cdot(\I_{n^2} - \BPhi\BPhi^{\top})
\]
and thus it remains to show that $(\BPi - \bar{\BP})(\I_{n^2} - \BPhi\BPhi^{\top})(\BPi - \bar{\BP})\succeq c'(\BPi - \bar{\BP})$ for some $c'>0$. Since
\begin{align*}
& (\BPi - \bar{\BP})(\I_{n^2} - \BPhi\BPhi^{\top})(\BPi - \bar{\BP}) = (\BPi - \bar{\BP}) - (\BPi - \bar{\BP}) \BPhi\BPhi^{\top}(\BPi - \bar{\BP}),
\end{align*}
it remains to show that 
$\lambda_{\max}( (\BPi - \bar{\BP})  \BPhi\BPhi^{\top} (\BPi - \bar{\BP}))\leq 1-c'$
and equivalently
\[
\lambda_{\max}( \BPhi^{\top}(\BPi - \bar{\BP})\BPhi ) \leq 1-c'\Longleftrightarrow
 \lambda_{\min} \left(\BPhi^{\top}\left(\frac{\bx\bx^{\top}}{n} + \bar{\BP}\right)\BPhi  \right) \geq c' 
\]
which follows from
\[
\BPhi^{\top}(\BPi - \bar{\BP})\BPhi = \I_n - \BPhi^{\top}\left(\frac{\bx\bx^{\top}}{n} + \bar{\BP}\right)\BPhi.
\]

Since $\bar{\BP} = n^{-2}(\BJ_n\otimes\I_n - \I_n\otimes\BJ_n)^2$, we first compute
\begin{align*}
(\BJ_n\otimes\I_n - \I_n\otimes\BJ_n) \BPhi & = (\BJ_n\otimes \I_n -  \I_n \otimes\BJ_n )[\bu_1\otimes\bu_1,\cdots,\bu_n\otimes\bu_n] \\
& = \left( [\bone_n\otimes\bu_1,\cdots,\bone_n\otimes\bu_n] - [\bu_1\otimes \bone_n,\cdots,\bu_n\otimes\bone_n]\right)\diag(\BU^{\top}\bone_n)\\
& = (\bone_n\otimes\BU - \BU\otimes\bone_n) \diag(\BU^{\top}\bone_n) 
\end{align*}
and then
\begin{align*}
\BPhi^{\top}\bar{\BP}\BPhi & = \frac{1}{n^2} \diag(\BU^{\top}\bone_n)   (\bone_n^{\top}\otimes\BU^{\top} - \BU^{\top}\otimes\bone_n^{\top})(\bone_n\otimes\BU - \BU\otimes\bone_n) \diag(\BU^{\top}\bone_n)  \\
& = \frac{1}{n^2} \diag(\BU^{\top}\bone_n)   (2n \I_n - \BU^{\top}\bone_n\otimes\bone_n^{\top}\BU - \bone_n^{\top}\BU\otimes\BU^{\top}\bone_n) \diag(\BU^{\top}\bone_n)  \\
& = \frac{2}{n} \diag(\BU^{\top}\bone_n)   ( \I_n - \BU^{\top}\BP\BU ) \diag(\BU^{\top}\bone_n).
\end{align*}
where $\BU^{\top}\bone_n \otimes \bone_n^{\top}\BU = \BU^{\top}\BJ_n\BU.$

Using~\eqref{eq:phiphix}, we have
\begin{align*}
\BPhi^{\top}(\bar{\BP} + n^{-1}\bx\bx^{\top})\BPhi & = \frac{2}{n} \diag(\BU^{\top}\bone_n)   ( \I_n - \BU^{\top}\BP\BU ) \diag(\BU^{\top}\bone_n)  + \frac{\BJ_n}{n}  
\end{align*}
where $( \I_n - \BU^{\top}\BP\BU ) \diag(\BU^{\top}\bone_n) \bone_n = ( \I_n - \BU^{\top}\BP\BU ) \BU^{\top}\bone_n = 0$. Therefore, we only need to control the second smallest eigenvalue of the first term.

Since $\BU^{\top}\BP\BU$ is rank-1 and $\Tr(\BU^{\top}\BP\BU)= 1$, we have $\I_n - \BU^{\top}\BP\BU$ is a projection matrix.
The second smallest eigenvalue of $\diag(\BU^{\top}\bone_n)   ( \I_n - \BU^{\top}\BP\BU ) \diag(\BU^{\top}\bone_n) $ is lower bounded by $\min|\lag \bu_i, \bone_n\rag|^2.$
Therefore, we have
\begin{equation}\label{def:cpp}
\lambda_{\min}(\BPhi^{\top}(\bar{\BP} + n^{-1}\bx\bx^{\top})\BPhi ) \geq  \min\left\{\frac{2 |\lag \bu_i,\bone_n\rag|^2}{n}, 1\right\} \geq \frac{2\min |\lag \bu_i,\bone_n\rag|^2}{n} = :c'.
\end{equation}
where $\min_i |\lag \bu_i,\bone_n\rag|^2 \leq 1$
and it holds
$(\BPi - \bar{\BP})(\I_{n^2} - \BPhi\BPhi^{\top})(\BPi - \bar{\BP})\succeq c'(\BPi - \bar{\BP}).$
Finally, we have
\begin{align*}
(\BPi - \bar{\BP}) \BM(\BPi - \bar{\BP}) & \succeq \min_{i\neq j}(\lambda_i - \lambda_j)^2 (\BPi - \bar{\BP}) (\I_{n^2} - \BPhi\BPhi^{\top})(\BPi - \bar{\BP}) \\
& \succeq c' \min_{i\neq j}(\lambda_i - \lambda_j)^2  (\BPi - \bar{\BP}).
\end{align*}
Lemma~\ref{lem:supp} implies that
\[
\BM+ nt\bar{\BP} \succeq \left( \frac{2}{n}\min |\lag \bu_i, \bone_n\rag|^2\cdot \min_{i\neq j}(\lambda_i - \lambda_j)^2 - \eps \right)(\I_{n^2} - \BPi)
\] 
for a sufficiently large $t$ where $\eps > 0$ can be arbitrarily small for a sufficiently large $t$.
\end{proof}

\begin{proof}[\bf Proof of Theorem~\ref{thm:main}]

Theorem~\ref{thm:suff} indicates that if the second smallest eigenvalue of $\BQ$ is nonnegative, then the exactness holds, i.e., $\BX = \bx\bx^{\top}$ is a global minimizer to~\eqref{def:sdr}. 

In presence of noise, $\BM$ equals
\begin{align*}
\BM & = (\I_n\otimes\BA - \BA\otimes\I_n - \BDelta\otimes\I_n)^2 \\
& = (\I_n\otimes\BA - \BA\otimes\I_n)^2-2 \BDelta\otimes\BA + ( \BDelta\BA +\BA\BDelta + \BDelta^2)\otimes\I_n.
\end{align*}
Let 
\begin{equation}\label{def:dualcert}
\begin{aligned}
\BT & =-t\BJ_n + \BDelta^2  + \BDelta\BA+\BA\BDelta-2\ddiag(\BA\BDelta)  -c(\BJ_n - \I_n), \\
\BK & =-t\BJ_n, ~~\BH = t\BJ_n/n, ~~\BZ  = t\BJ_n/n, 
\end{aligned}
\end{equation}
for some $c$ and $t > 0$. Now we verify~\eqref{eq:Bxdiag} in Proposition~\ref{prop:consBQ}:
\begin{align*}
\diag(\BT+\BK) + (\BZ+\BH)\bone_n & = -t\bone_n + \diag(\BDelta^2) - t\bone_n + 2t\bone_n = \diag(\BDelta^2)
\end{align*}
where $\diag(\BT) = -t\bone_n + \diag(\BDelta^2).$ Therefore, we can choose $\BB$ in the form of~\eqref{def:B}. To ensure the KKT condition~\eqref{eq:kktfinal_simp}, it suffices to have~\eqref{def:MSx} holds so that $\BB\geq 0$ and $\BQ$ in~\eqref{def:Qcons} is positive semidefintie. 
For the first requirement, we have
\begin{align*}
\BT + \BJ_n \BZ + \BK +\BH\BJ_n & = -c(\BJ_n - \I_n)+ \BDelta^2  + \BDelta\BA+\BA\BDelta-2\ddiag(\BA\BDelta)
\end{align*}
and thus $\BB\geq 0$ is guaranteed if
\[
\BDelta^2 +\BDelta\BA-\BA\BDelta - ( \BT + \BJ_n\BZ + \BK + \BH\BJ_n) = c (\BJ_n - \I_n) + 2\ddiag(\BA\BDelta) - 2\BA\BDelta \geq 0
\]
provided that $c + 2 \min_{i\neq j}\lag \BA_i, \BDelta_j\rag \geq 0$. 

For $\BQ\succeq 0$, we note 
{\small
\begin{align*}
\BS
& = (\BDelta^2 + \BDelta\BA+\BA\BDelta-2\ddiag(\BA\BDelta)   - c(\BJ_n - \I_n) )\otimes\I_n -  nt\bar{\BP}
\end{align*}}
where $\bar{\BP} : = n^{-2}(\BJ_n\otimes\I_n - \I_n \otimes\BJ_n )^2.$
Then 
\begin{align*}
 \BM - \BS & = (\I_n\otimes\BA - \BA\otimes\I_n)^2-2 \BDelta\otimes\BA + ( \BDelta\BA +\BA\BDelta + \BDelta^2)\otimes\I_n \\
& \qquad - (\BDelta^2 + \BDelta\BA+\BA\BDelta-2\ddiag(\BA\BDelta)   - c(\BJ_n - \I_n) )\otimes\I_n +  nt\bar{\BP} \\
& = (\I_n\otimes\BA - \BA\otimes\I_n)^2 +  nt\bar{\BP}  - 2\BDelta\otimes \BA + (2\ddiag(\BA\BDelta)   + c(\BJ_n - \I_n) )\otimes \I_n.
\end{align*}

By choosing $\BB$ in the form of~\eqref{def:B}, we have $\BQ = (\I_{n^2} - n^{-1}\bx\bx^{\top})(\BM - \BS)(\I_{n^2} - n^{-1}\bx\bx^{\top})$ and the second smallest eigenvalue satisfies 
\begin{equation}\label{eq:lambda2Q}
\begin{aligned}
\lambda_{2}\left(\BQ\right) & \geq \lambda_2((\I_n\otimes\BA - \BA\otimes\I_n)^2 +  nt\bar{\BP}) - 2\|\BDelta\otimes\BA\| + 2 \min_{1\leq i\leq n} \lag \BA_i, \BDelta_i\rag -c \\
\text{(Theorem~\ref{thm:noisefree})}~~& \geq  \frac{2}{n}\min_{i\neq j}(\lambda_i - \lambda_j)^2 \min |\lag \bu_i,\bone_n\rag|^2 - 2\|\BDelta\| \|\BA\| - 2 \|\BA\BDelta\|_{\max}
\end{aligned}
\end{equation}
where $\|\BDelta\otimes\BA\| = \|\BA\|\|\BDelta\|$, $c=-2\min_{i\neq j}\lag \BA_i, \BDelta_j\rag$ is chosen, and $t$ is sufficiently large.
Then by Theorem~\ref{thm:suff}, the exactness holds if $\lambda_2(\BQ) > 0$, i.e., $\BX=\bx\bx^{\top}$ is the unique global minimizer to~\eqref{def:sdr}.
\end{proof}

\section{Proof of Theorem~\ref{thm:main2}}\label{s:proof2}

\subsection{Dual program and optimality condition of~\eqref{def:sdr2}}
We start with deriving the dual form of~\eqref{def:sdr2}.
For each constraint, we assign a dual variable:
\begin{align*}
\begin{bmatrix}
\BQ & \bq \\
\bq & z
\end{bmatrix}\succeq 0 &: \begin{bmatrix}
\BX & \diag(\BX) \\
\diag(\BX)^{\top} & 1 
\end{bmatrix}\succeq 0, \\
\BB \geq 0 &: \BX \geq 0, \\
\BT\in\RR^{n\times n} &: \Tr(\BX_{ij}) = \delta_{ij}, ~~1\leq i,j\leq n, \\
\BZ\in\RR^{n\times n} &: \lag \BX_{ij}, \BJ_{n}\rag = 1, ~~1\leq i,j\leq n,\\
\BK\in\RR^{n\times n}&: \sum_{i=1}^n \BX_{ii} = \I_n, \\
\BH\in \RR^{n\times n} &: \sum_{i,j} \BX_{ij} = \BJ_n, \\
\bmu\in\RR^n &: \mat(\diag(\BX))\bone_n = \bone_n, \\
\blambda\in \RR^n &:  \mat(\diag(\BX))^{\top}\bone_n = \bone_n,
\end{align*} 
where $\diag(\BX)$ takes the diagonal entries of $\BX$ and forms them into a column vector. 

For the constraint $\mat(\diag(\BX))$, 
\begin{align*}
\bmu^{\top}(\mat(\diag(\BX)) \bone_n - \bone_n) 
& = \lag \mat(\diag(\BX)), \bmu\bone_n^{\top}\rag - \lag \bmu,\bone_n\rag\\
& =\lag \diag(\BX), \bone\otimes \bmu \rag - \lag \bmu,\bone_n\rag = \lag \I_n\otimes \diag( \bmu), \BX\rag -\lag \bmu, \bone_n\rag \\
\blambda^{\top}(\mat(\diag(\BX))^{\top} \bone_n - \bone_n) 
& = \lag \mat(\diag(\BX)), \bone_n\blambda^{\top}\rag - \lag \blambda,\bone_n\rag\\
& =\lag \diag(\BX), \blambda\otimes \bone_n \rag - \lag \blambda,\bone_n\rag = \lag \diag( \blambda)\otimes \I_n, \BX\rag -\lag \blambda, \bone_n\rag 
\end{align*} 
where $\VEC(\bu\bv^{\top}) = \bv\otimes \bu.$ Also
\[
\left\lag 
\begin{bmatrix}
\BX & \diag(\BX) \\
\diag(\BX)^{\top} & 1 
\end{bmatrix},
\begin{bmatrix}
\BQ & \bq \\
\bq & z
\end{bmatrix}
\right\rag = \lag \BX,\BQ\rag + 2\lag\BX, \diag(\bq)\rag + z 
\]
where $\lag\diag(\BX),\bq \rag = \lag\BX,\diag(\bq)\rag.$

Now the Lagrangian function is {\small
\begin{align*}
& {\cal L}(\BX,\BQ,\BB,\BT,\BZ,\BK,\BH,\bmu,\blambda,\bq,z) \\
& \qquad : = \lag \BX,\BM-\BQ- 2\diag(\bq)-\BB  -( \I_n\otimes \diag(\bmu) + \diag(\blambda)\otimes \I_n) \rag  + \lag \blambda+\bmu,\bone_n\rag\\
& \qquad \quad - \sum_{i,j} t_{ij}(\Tr(\BX_{ij}) -\delta_{ij})  - \sum_{i,j} z_{ij}(\lag \BX_{ij},\BJ_n\rag - 1)   -\left \lag \sum_{i=1}^n \BX_{ii} - \I_n, \BK\right\rag - \left\lag \sum_{i,j}\BX_{ij} - \BJ_n, \BH\right\rag -z \\
&  \qquad = \left\lag  \BX,\BM-\BQ-2\diag(\bq)-\BB-\diag(\bone\otimes \bmu + \blambda\otimes \bone_n)  - \BT\otimes \I_n - \BZ\otimes \BJ_n -  \I_n\otimes \BK -  \BJ_n\otimes \BH\right\rag \\
&\qquad\quad   + \Tr(\BT + \BK) + \lag \BJ_n,\BZ+\BH\rag + \lag \blambda+\bmu\rag - z.
\end{align*}}
We define
\begin{equation}\label{def:Lambda}
\BLambda = \I_n \otimes\diag(\bmu) + \diag(\blambda)\otimes \I_n, ~~~~
\BS =  \BT\otimes \I_n +  \I_n\otimes \BK + \BJ_n\otimes \BH+ \BZ\otimes \BJ_n.
\end{equation}
Then the Lagrangian equals
\[
{\cal L}(\cdot) = \lag \BX, \BM - \BQ-2\diag(\bq) -\BB - \BLambda - \BS \rag + \Tr(\BT + \BK) + \lag \BJ_n,\BZ+\BH\rag + \lag \blambda+\bmu, \bone_n\rag - z.
\] 
The resulting dual program becomes
\begin{equation}
\begin{aligned}
\max~ \qquad &\Tr(\BT + \BK) + \lag \BJ_n,\BZ+\BH\rag + \lag \blambda+\bmu, \bone_n\rag - z,  \\
\text{s.t.} ~\qquad
& \BM=\BQ+2\diag(\bq) +\BB+ \BLambda + \BS,  \\
& \begin{bmatrix}
\BQ & \bq \\
\bq & z
\end{bmatrix}\succeq 0,  ~~\BB\geq 0. 
\end{aligned} 
\end{equation}
As a result, the KKT conditions are
\begin{enumerate}
\item Stationarity:
\[
\BM=\BQ+ 2\diag(\bq)+ \BB  + \BLambda + \BS,
\]

\item Dual feasibility:
\[
\begin{bmatrix}
\BQ & \bq \\
\bq^{\top} & z
\end{bmatrix} \succeq 0,~~~~\BB \geq 0,
\]

\item Complementary slackness:
\[
\left\lag 
\begin{bmatrix}
\BX & \diag(\BX) \\
\diag(\BX)^{\top} & 1 
\end{bmatrix},
\begin{bmatrix}
\BQ & \bq \\
\bq & z
\end{bmatrix}
\right\rag = 0,~~~\lag \BB, \BX\rag = 0.
\]
\end{enumerate}

Suppose $\BX = \bx\bx^{\top}$ is a global minimizer, and then we can try to simplify the KKT optimality condition by determining some of the dual variables. Using the complementary slackness condition gives rise to 
\begin{align*}
\begin{bmatrix}
\BQ & \bq \\
\bq^{\top} & z
\end{bmatrix} \begin{bmatrix}
\bx \\
1
\end{bmatrix} =0,~~~~\lag \BB, \bx\bx^{\top}\rag  = 0,
\end{align*}
which leads to $\BQ\bx + \bq  = 0,z= -\lag \bq,\bx\rag,$ and $\BB\circ \bx\bx^{\top} =0.$
Now we look into the KKT conditions again:
for the dual feasibility, it holds that
\[
\begin{bmatrix}
\BQ & -\BQ\bx \\
-(\BQ\bx)^{\top} & \bx^{\top}\BQ\bx 
\end{bmatrix} \succeq 0  ~ \Longleftrightarrow~
 \BQ\succeq 0.
\]
For the stationarity, we have
\[
\BM = \BQ - 2\diag(\BQ\bx) + \BB + \BLambda + \BS.
\]
This  linear system $\BQ - 2\diag(\BQ\bx) = \BM - (\BB + \BLambda +\BS)$ has a unique solution for $\BQ$:
\[
\BQ\bx - 2\diag(\bx)\BQ\bx = (\BM - \BB  -\BLambda -\BS)\bx \Longrightarrow \BQ\bx = (\I - 2\diag(\bx))(\BM - \BB  - \BLambda -\BS)\bx
\]
and then
\begin{align*}
\BQ & =  \BM - (\BB + \BLambda +\BS) + 2\diag(\BQ\bx) \\
& =  \BM - (\BB + \BLambda +\BS) + 2(\I - 2\diag(\bx))\diag((\BM - \BB  - \BLambda -\BS)\bx)
\end{align*}
Therefore, the KKT condition becomes
\begin{equation}\label{eq:kktfull2}
\begin{aligned}
& \BQ = \BM - \BB  - \BLambda -\BS + 2(\I - 2\diag(\bx))\diag((\BM - \BB  - \BLambda -\BS)\bx)\succeq 0, \\
& \BB \geq 0,~\BB\circ \bx \bx^{\top} = 0,
\end{aligned}
\end{equation}
where $\BS$ is defined in~\eqref{def:Lambda}. 
Next, we show that~\eqref{eq:kktfull2} implies the (unique) global optimality of $\BX.$
\begin{theorem}\label{thm:suff2}
Suppose that there exist $\BB$, $\BS$, and $\BLambda$ such that~\eqref{eq:kktfull2} holds, then $\BX= \bx\bx^{\top}$ is a global minimizer. Moreover, if $\BQ \succ 0$, then $\BX$ is the unique global minimizer.
\end{theorem}

\begin{proof}[\bf Proof of Theorem~\ref{thm:suff2}]
Suppose that
\[
\BQ = \BM - \BB - \BLambda - \BS +  2\diag(\BQ\bx) \succeq 0
\]
which gives
\[
\bx^{\top} \BQ\bx = -\bx^{\top}( \BM - \BB - \BLambda - \BS)\bx.
\]
For any feasible solution $\widehat{\BX}$, we have
\begin{equation}\label{eq:csopt}
\begin{aligned}
& \left\lag 
\begin{bmatrix}
\BQ & - \BQ\bx \\
- (\BQ\bx)^{\top} & \bx^{\top}\BQ\bx 
\end{bmatrix}, 
\begin{bmatrix}
\widehat{\BX} & \diag(\widehat{\BX}) \\
\diag(\widehat{\BX})^{\top} & 1
\end{bmatrix}
\right\rag 
 = \lag \widehat{\BX}, \BQ -2\diag(\BQ\bx)\rag + \bx^{\top}\BQ\bx \\
& \qquad = \lag \widehat{\BX} - \bx\bx^{\top}, \BM - \BB-\BLambda- \BS\rag  = \lag \widehat{\BX} - \bx\bx^{\top}, \BM\rag - \lag \widehat{\BX}, \BB\rag \geq 0
\end{aligned}
\end{equation}
where $\lag \widehat{\BX} - \bx\bx^{\top}, \BLambda+ \BS\rag = 0$ follows from the feasibility of $\widehat{\BX}$ and $\bx\bx^{\top}$, and $\lag \bx\bx^{\top},\BB\rag = 0.$ Therefore, it holds $\lag \widehat{\BX} - \bx\bx^{\top}, \BM\rag \geq \lag \widehat{\BX}, \BB\rag\geq 0$
which implies $\BX = \bx\bx^{\top}$ is a global minimizer to~\eqref{def:sdr2}. In particular, if $\BQ\succ 0$, then $\BX$ is the unique global minimizer since~\eqref{eq:csopt} becomes strictly positive for any $\widehat{\BX}\neq \bx\bx^{\top}.$
\end{proof}

\begin{proof}[\bf Proof of Theorem~\ref{thm:main2}]
Consider the dual certificate $\BT,\BK,\BZ$ and $\BH$ in~\eqref{def:dualcert}, and then it holds that
\[
\diag(\BT+\BK) + (\BZ+\BH)\bone_n = \diag(\BDelta^2),~~~\BB\bx = (\BM-\BS)\bx
\]
where $\BB$ is chosen in the form of~\eqref{def:B}.
As a result,  $\BQ$ satisfies
\begin{align*}
\BQ & = \BM - \BB  - \BLambda -\BS + 2(\I - 2\diag(\bx))\diag((\BM - \BB  - \BLambda -\BS)\bx) \\
& =\BPi(\BM   -\BS)\BPi - \left( \BLambda + 2(\I - 2\diag(\bx))\diag( \BLambda \bx) \right) 
\end{align*}
where $\BB\bx = (\BM-\BS)\bx$ and $\BPi := \I - \bx\bx^{\top}/n.$
Note that~\eqref{eq:lambda2Q} gives
\[
\lambda_2\left( \BPi(\BM   -\BS)\BPi \right) \geq \frac{2}{n}\min_{i\neq j}(\lambda_i - \lambda_j)^2  \min |\lag \bu_i,\bone_n\rag|^2 - 2\|\BDelta\| \|\BA\| -2 \|\BA\BDelta\|_{\max}.
\]
Note that $\mat(\BLambda\bx) = \diag(\blambda+\bmu)$ where $\BLambda$ is defined in~\eqref{def:Lambda}.
Also we have
\[
\mat(\diag(\BLambda) + 2(\I - 2\diag(\bx)) \BLambda \bx) = \bone_n\blambda^{\top} + \bmu\bone_n^{\top} - 2\diag(\blambda+\bmu).
\]
In particular, if we choose $\blambda$ and $\bmu$ as constant vectors, i.e., $\blambda = \bmu = t\bone_n$, then
\[
\mat(\diag(\BLambda) + 2(\I - 2\diag(\bx)) \BLambda \bx)  = 2t( \BJ_n - 2\I_n )
\]
and thus 
\begin{align*}
\BQ  & = \BPi(\BM   -\BS)\BPi - 2t (\I - \diag(\bx)) + 2t\diag(\bx) \\
& \succeq \lambda_2(\BPi(\BM   -\BS)\BPi) \left(\I - \frac{\bx\bx^{\top}}{n}\right)- 2t (\I - \diag(\bx)) + 2t\diag(\bx)
\end{align*}
For any $0 < 2t< \lambda_2(\BPi(\BM-\BS)\BPi)$, then $\BQ\succ 0$ holds and therefore $\BX=\bx\bx^{\top}$ is the unique global minimizer to~\eqref{def:sdr2}, following from Theorem~\ref{thm:suff2}.
\end{proof}

\section*{Appendix}

\begin{lemma}\label{lem:supp}
Let $\BPi$ and $\BV$ be two matrices of same size and $\BP_{\Pi}$ and $\BP_V$ be their orthogonal projection matrices respectively that satisfy $\BP_V(\I - \BP_{\Pi})= 0$, i.e., $\Ran(\BV) \subseteq \Ran(\BPi)$. Suppose
\[
(\BP_{\Pi}-\BP_V)\BA (\BP_{\Pi}-\BP_V)\succeq c(\BP_{\Pi}-\BP_V)
\]
for some $c > 0$. Then for any $0 < c'' < c$, we have 
\[
\BP_{\Pi}\BA \BP_{\Pi} + t\BP_V\succeq c''\BP_{\Pi}
\] 
for sufficiently large $t$.

\end{lemma}
\begin{proof}[\bf Proof of Lemma~\ref{lem:supp}]
We decompose $\BP_{\Pi}\BA\BP_{\Pi}+t\BP_V$ into
\[
\BP_{\Pi}\BA\BP_{\Pi}+t\BP_V = (\BP_{\Pi}-\BP_V)\BA (\BP_{\Pi}-\BP_V) + (\BP_{\Pi} - \BP_V)\BA\BP_V + \BP_V \BA(\BP_{\Pi} - \BP_V) + (\BP_V\BA\BP_V + t\BP_V).
\]
Let $\BR := \BP_{\Pi} - (\BP_{\Pi} - \BP_V)\BA (\BP_V \BA\BP_V + t\BP_V)^{\dagger}$
where $\Ran(\BR) \subseteq \Ran(\BPi)$, ``$\dagger$" stands for the Moore-Penrose pseudo-inverse and
\begin{equation}\label{def:BPV}
(\BP_V \BA\BP_V + t\BP_V)(\BP_V \BA\BP_V + t\BP_V)^{\dagger} = \BP_V.
\end{equation}
Note that
\[
(\BR - \BP_{\Pi})^2 = (\BP_{\Pi} - \BP_V)\BA (\BP_V \BA\BP_V + t\BP_V)^{\dagger} (\BP_{\Pi} - \BP_V)\BA (\BP_V \BA\BP_V + t\BP_V)^{\dagger} = 0
\]
where $(\BP_V \BA\BP_V + t\BP_V)^{\dagger} (\BP_{\Pi} - \BP_V) = 0$. Using $\BP_{\Pi}\BR = \BR\BP_{\Pi}= \BR$ and $(\BR - \BP_{\Pi})^2 = 0$ gives
\[
\BR^2 - 2\BR + \BP_{\Pi} = 0
\]
which implies $\BR$ is invertible on $\Ran(\BPi)$ as 1 is an eigenvalue of $\BR$ with multiplicity equal to the rank of $\BP_{\Pi}.$ Therefore, $\BP_{\Pi}\BA\BP_{\Pi}+t\BP_V\succ_{\Pi} 0 $ is equivalent to $\BR(\BP_{\Pi}\BA\BP_{\Pi}+t\BP_V)\BR^{\top}\succ_{\Pi} 0$.
Now 
\begin{align*}
& \BR(\BP_{\Pi}\BA\BP_{\Pi} + t\BP_V) \\
& = (\BP_{\Pi} - (\BP_{\Pi} - \BP_V)\BA (\BP_V \BA\BP_V + t\BP_V)^{\dagger})(\BP_{\Pi}\BA\BP_{\Pi} + t\BP_V) \\
& = \BP_{\Pi}\BA\BP_{\Pi} + t\BP_V - (\BP_{\Pi} - \BP_V)\BA (\BP_V \BA\BP_V + t\BP_V)^{\dagger} (\BP_V \BA(\BP_{\Pi} - \BP_V) + (\BP_V\BA\BP_V + t\BP_V)) \\
& = \BP_{\Pi}\BA\BP_{\Pi} + t\BP_V -  (\BP_{\Pi} - \BP_V)\BA (\BP_V \BA\BP_V + t\BP_V)^{\dagger} \BA (\BP_{\Pi} - \BP_V) - (\BP_{\Pi} - \BP_V)\BA\BP_V \\
& =  (\BP_{\Pi}-\BP_V)\BA (\BP_{\Pi}-\BP_V)  + \BP_V \BA(\BP_{\Pi} - \BP_V) + (\BP_V\BA\BP_V + t\BP_V) \\
& \qquad  -  (\BP_{\Pi} - \BP_V)\BA (\BP_V \BA\BP_V + t\BP_V)^{\dagger} \BA (\BP_{\Pi} - \BP_V)
\end{align*}
which follows from $\BP_{\Pi}\BP_{V} = \BP_V$ and~\eqref{def:BPV}.
Then
\begin{align*}
& \BR(\BP_{\Pi}\BA\BP_{\Pi} + t\BP_V) \BR^{\top} \\
& \qquad = (\BP_{\Pi}-\BP_V)\BA (\BP_{\Pi}-\BP_V)  + \BP_V \BA(\BP_{\Pi} - \BP_V) + (\BP_V\BA\BP_V + t\BP_V) \\
& \qquad ~~~ -  (\BP_{\Pi} - \BP_V)\BA (\BP_V \BA\BP_V + t\BP_V)^{\dagger} \BA (\BP_{\Pi} - \BP_V) - \BP_V\BA (\BP_{\Pi} - \BP_V)\\
&\qquad = (\BP_V\BA\BP_V + t\BP_V) + (\BP_{\Pi}-\BP_V)\BA (\BP_{\Pi}-\BP_V)    -  (\BP_{\Pi} - \BP_V)\BA (\BP_V \BA\BP_V + t\BP_V)^{\dagger} \BA (\BP_{\Pi} - \BP_V)
\end{align*}
where $(\BP_V\BA\BP_V + t\BP_V) (\BP_V \BA\BP_V + t\BP_V)^{\dagger} = \BP_V$ for sufficiently large $t$. 

Next, we will show that $\BR(\BP_{\Pi}\BA\BP_{\Pi} + t\BP_V) \BR^{\top}$ is strictly positive semidefinite when restricted to $\BP_V.$
Suppose $(\BP_{\Pi}-\BP_V)\BA (\BP_{\Pi}-\BP_V) \succeq c(\BP_{\Pi} - \BP_V)$, and then for a sufficiently large $t>0$, it holds for $0 <c'<c$ that  
\begin{align*} 
& (\BP_{\Pi}-\BP_V)\BA (\BP_{\Pi}-\BP_V)   -  (\BP_{\Pi} - \BP_V)\BA (\BP_V \BA\BP_V + t\BP_V)^{\dagger} \BA (\BP_{\Pi} - \BP_V) \succeq c' (\BP_{\Pi} - \BP_V)
\end{align*}
as the second term can be arbitrarily small for a sufficiently large $t>0.$
 Then
\[
\BR(\BP_{\Pi}\BA\BP_{\Pi} + t\BP_V) \BR^{\top} \succeq  c' \BP_V + c'(\BP_{\Pi} - \BP_V) \succeq c'\BP_{\Pi}
\]
where $\BP_V\BA\BP_V + t\BP_V\succeq c' \BP_{V}$ for  sufficiently large $t$.

Note that $\BR$ is invertible on $\Ran(\BPi)$ and also the range of $\BR$ belongs to $\Ran(\BPi).$ Therefore, we have
\[
\BP_{\Pi}\BA\BP_{\Pi} + t\BP_V \succeq c' (\BR \BR^{\top})^{\dagger} \succeq c'' \BP_{\Pi}
\]
for any $c'' < c$
where $\| \BR \BR^{\top} - \BP_{\Pi}\|$ goes to 0 at the rate of $1/t$ and thus $\BR \BR^{\top}$ is close to $\BP_{\Pi}$ for sufficiently large  $t$. 
\end{proof}

\bibliographystyle{abbrv}

\end{document}